\documentclass[11pt]{amsart}

\usepackage{amsmath}
\usepackage{amssymb}
\usepackage{amsthm}
\usepackage{geometry}
\usepackage{xcolor}
\definecolor{babyblue}{RGB}{13,139,74}
\usepackage{tikz}
\usepackage{comment}
\usepackage{float}

\usepackage[
    backend=biber,
    style=numeric,
    sorting=nyt,
    giveninits=true,
    maxbibnames=99,
    doi=false,
    eprint=false,
    isbn=false,
    url=false
]{biblatex}
\DeclareFieldFormat[article,incollection,inproceedings,online]{title}{#1\isdot}
\DeclareFieldFormat[article]{volume}{\mkbibbold{#1}}
\DeclareFieldFormat[article]{pages}{#1}
\renewbibmacro*{in:}{%
  \ifentrytype{article}
    {}
    {\printtext{\bibstring{in}\intitlepunct}}}
\renewbibmacro*{journal+issuetitle}{%
  \usebibmacro{journal}%
  \setunit*{\addcomma\space}%
  \iffieldundef{series}
    {}
    {\printfield{series}%
     \setunit{\addcomma\space}}%
  \printfield{volume}%
  \setunit{\addspace}%
  \usebibmacro{issue+date}%
  \setunit*{\addcolon\space}%
  \usebibmacro{issue}%
  \newunit}
\usepackage[
    colorlinks=true,
    citecolor=babyblue,
    linkcolor=babyblue,
    urlcolor=babyblue,
]{hyperref}

\usepackage{pgfplots}
\pgfplotsset{compat=1.18}

\usetikzlibrary{arrows.meta}
\theoremstyle{plain}
\newtheorem{theorem}{Theorem}[section]
\newtheorem{lemma}[theorem]{Lemma}
\newtheorem{proposition}[theorem]{Proposition}
\newtheorem{corollary}[theorem]{Corollary}
\newtheorem{conjecture}[theorem]{Conjecture}
\newtheorem{fact}[theorem]{Fact}

\newtheorem{definition}[theorem]{Definition}

\newtheorem{remark}[theorem]{Remark}

\newcommand{\R}{\mathbb R}

\newcommand{\N}{\mathbb N}
\newcommand{\Z}{\mathbb Z}
\newcommand{\Q}{\mathbb Q}
\newcommand{\mG}{\mathcal G}
\renewcommand{\P}{\mathbb P}
\newcommand{\E}{\mathbb E}

\title[On Structures of Nearest-Predecessor Random Graphs]
{Local Limits and Structures of \\Nearest-Predecessor Random Graphs}

\author{Xiaoyao Huang}
\address{Department of Mathematics, University of Michigan, Ann Arbor, MI 48109}
\email{xyrushac@umich.edu}

\subjclass[2020]{Primary 05C80; Secondary 60C05, 05C63, 60F05}
\keywords{Random graphs, Benjamini--Schramm convergence, degree, one-ended, recurrence}

\begin{document}
\begin{abstract}
We study the two-sided nearest-predecessor graph on \([n]\), where vertices arrive in a uniformly random order and each new vertex connects to its nearest earlier vertex on each side.
We construct an infinite graph on \(\Z\) from i.i.d.\ continuous insertion times and prove quenched Benjamini--Schramm convergence to it almost surely. The limit is almost surely one-ended and recurrent.
We resolve an open problem of Ponomarenko by proving that the maximum degree is asymptotic in probability to \(\log n/\log(3/2)\) as \(n\to\infty\).
We obtain the probability generating function and a CLT for the number of edges crossing the boundary of \([n]\) in the infinite limit.

We also study the \(K\)-nearest-predecessor graph, in which each new vertex attaches to its \(K\) nearest earlier vertices without regard to side.
For every fixed \(K\geq1\), we construct the infinite limit, prove almost sure quenched Benjamini--Schramm convergence, and show that it is almost surely connected, one-ended, and recurrent, as well as nonplanar when \(K\geq2\).
\end{abstract}
\maketitle


\section{Introduction}
For \(n\geq1\), let \([n]=\{1,\ldots,n\}\).
We construct random graphs on the ordered vertex set \([n]\) from a uniformly random permutation \(\tau_n\) of \([n]\).
The value \(\tau_n(v)\) is the \emph{insertion time} of \(v\); thus vertices are inserted in increasing order of their \(\tau_n\) values.
A vertex \(u\) is a \emph{predecessor} of \(v\) if \(\tau_n(u)<\tau_n(v)\).
We study two graph models generated by this random order: the \emph{two-sided nearest-predecessor graph}, in which each new vertex connects to its nearest predecessor on the left and on the right, when they exist, and the \emph{\(K\)-nearest-predecessor graph}, in which each new vertex connects to its \(K\) nearest predecessors without a prescribed left-right allocation.

After reparametrizing the insertion times, the two-sided model is exactly the \emph{horizontal visibility graph} (HVG) of an i.i.d.\ sequence.
Luque, Lacasa, Ballesteros, and Luque~\cite{LuqueLacasaBallesterosLuque2010} derived the degree law \(p_d=\frac13(\frac23)^{d-2}\) for degree \(d\geq2\) and i.i.d.\ time series.
Gutin, Mansour, and Severini~\cite{GutinMansourSeverini2011} proved that a graph is an HVG if and only if it is outerplanar and has a Hamilton path.
Ponomarenko~\cite{Ponomarenko2026GreedyRouting} used record statistics to prove logarithmic greedy routing from \(1\) to \(n\) with high probability.
Beyond these results, sharper asymptotic structural analysis remains desirable.
In particular, maximum degree asymptotics were left open in the latter work; we solve this problem here (Theorem~\ref{thm:max-degree-asymptotics}).

The \(K\)-nearest-predecessor extension was suggested in~\cite{Ponomarenko2026GreedyRouting}; here we formalize and study the model for the first time.
For \(K=1\), the \(K\)-nearest-predecessor construction is closely related to the random on-line nearest-neighbor graph~\cite{PenroseWade2008}.

Our main contribution is a probabilistic framework, based on infinite local limits, for analyzing the structure and geometry of these models.
Quenched Benjamini--Schramm convergence provides a rigorous bridge between the finite and infinite graphs, allowing the neighborhood of a uniformly chosen finite vertex to be understood through its infinite limit.
The limiting graphs also reveal interesting features that arise only in the infinite setting.

We highlight four main results of this paper.
For the two-sided model, Theorem~\ref{thm:max-degree-asymptotics} proves that the maximum degree is asymptotic in probability to \(\log n/\log(3/2)\), while Theorem~\ref{thm:boundary-fluctuations} gives the generating function and a CLT for the number of edges crossing the boundary of \([n]\).
Theorems~\ref{thm:quenched-BS} and~\ref{thm:directionless-quenched-BS} establish almost sure quenched Benjamini--Schramm convergence for the two models.

\subsection{Two-Sided Nearest-Predecessor Graph}
\begin{definition}
\label{def:finite-two-sided-nearest-predecessor}
The graph \(G_n=([n],E_n)\) is obtained by adjoining an edge from each vertex \(v\) to its nearest predecessor \(u<v\), if one exists, and an edge from \(v\) to its nearest predecessor \(u>v\), if one exists.
\end{definition}

\begin{figure}[ht]
\centering
\begin{tikzpicture}[x=0.52cm,y=-0.35cm,line cap=round,line join=round]
\node[anchor=south west,font=\scriptsize] at (1,1) {insertion time};
\node[anchor=south,font=\small] at (9,0.5) {two-sided nearest-predecessor graph};

\draw (0,1) -- (0,17);
\foreach \t in {1,...,17} {
    \draw (-0.04,\t) -- (0.04,\t);
    \node[anchor=east,font=\scriptsize] at (-0.10,\t) {\t};
}

\draw (1,18) -- (17,18);
\foreach \v in {1,...,17} {
    \draw (\v,17.9) -- (\v,18.1);
    \node[anchor=base,font=\scriptsize] at (\v,19) {\v};
}
\node[anchor=base,font=\scriptsize] at (9,20) {vertex number \(v\)};

\foreach \v/\t in {1/7,2/8,3/11,4/5,5/14,6/13,7/1,8/6,9/10,10/4,11/2,12/9,13/15,14/16,15/3,16/12,17/17} {
    \node[anchor=base,font=\scriptsize] at (\v,20.9) {\t};
}
\node[anchor=base,font=\scriptsize] at (9,21.9) {insertion permutation \(\tau_n(v)\)};

\foreach \v in {1,...,17} {
    \draw[babyblue!25,line width=0.8pt] (\v,1) -- (\v,17);
}
\foreach \t in {1,...,17} {
    \draw[babyblue!25,line width=0.8pt] (1,\t) -- (17,\t);
}

\foreach \xa/\ya/\xb/\yb in {
1/7/2/8,
1/7/4/5,
2/8/3/11,
2/8/4/5,
3/11/4/5,
4/5/5/14,
4/5/6/13,
4/5/7/1,
5/14/6/13,
6/13/7/1,
7/1/8/6,
7/1/10/4,
7/1/11/2,
8/6/9/10,
8/6/10/4,
9/10/10/4,
10/4/11/2,
11/2/12/9,
11/2/15/3,
12/9/13/15,
12/9/15/3,
13/15/14/16,
13/15/15/3,
14/16/15/3,
15/3/16/12,
16/12/17/17
} {
    \draw[line width=1.05pt] (\xa,\ya) -- (\xb,\yb);
}

\foreach \v/\t in {1/7,2/8,3/11,4/5,5/14,6/13,7/1,8/6,9/10,10/4,11/2,12/9,13/15,14/16,15/3,16/12,17/17} {
    \node[draw,circle,fill=white,inner sep=1.2pt,font=\scriptsize] at (\v,\t) {\v};
}
\end{tikzpicture}
\caption{A drawing of a realization of the finite two-sided nearest-predecessor graph for \(n=17\), using the permutation
\((\tau_n(1),\ldots,\tau_n(17)) = (7,8,11,5,14,13,1,6,10,4,2,9,15,16,3,12,17)\).}
\label{fig:finite-graph-example}
\end{figure}
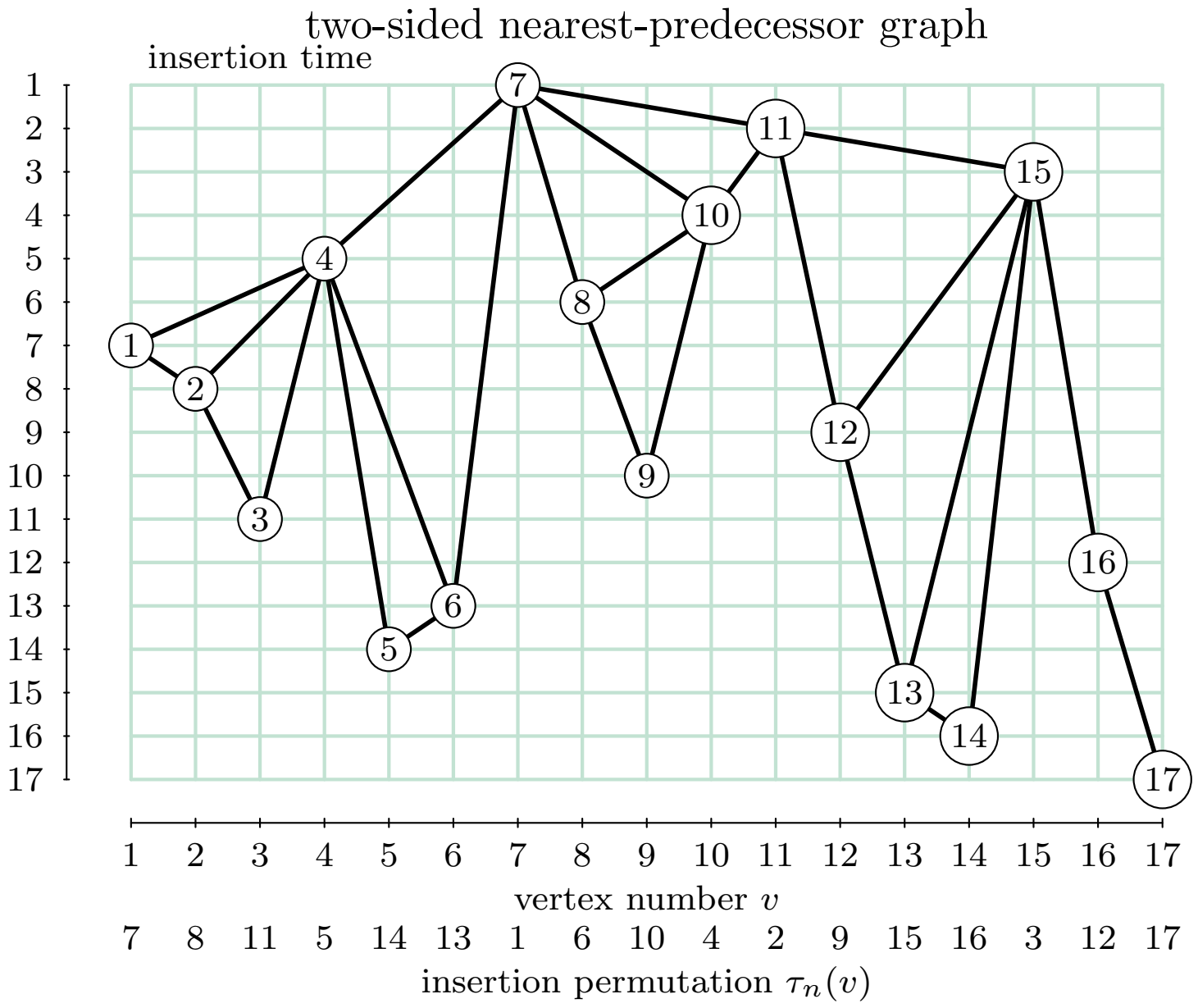

See Figure~\ref{fig:finite-graph-example} for an example.

We study the limiting local neighborhood of a uniformly chosen vertex in the Benjamini--Schramm sense~\cite{BenjaminiSchramm2001,AldousLyons2007}.

To describe the limiting graph, let \((T_i)_{i\in\Z}\) be i.i.d.\ random variables with distribution \(\operatorname{Unif}(0,1)\), where a smaller value represents an earlier insertion.
We work on the event of probability one on which all insertion times are distinct.

\begin{definition}
\label{def:infinite-two-sided-nearest-predecessor}
Define \(G_\infty=(\Z,E_\infty)\) as follows: for \(i<j\),
\begin{equation}
\label{inf-rule}
\{i,j\}\in E_\infty
\quad\Longleftrightarrow\quad
\max\{T_i,T_j\}<\min_{i<k<j}T_k.
\end{equation}
Here, the minimum over an empty set is \(+\infty\).
\end{definition}

The subgraph of \(G_\infty\) induced on \([n]\) is identically distributed with \(G_n\).
Our first goal is to show that \(G_\infty\) describes the limiting local neighborhood of a uniformly chosen vertex of \(G_n\).

For a rooted graph \((G,o)\), let \(B_r(G,o)\) denote the rooted subgraph induced by the vertices at graph distance at most \(r\) from \(o\).

Since the law of \(G_\infty\) is invariant under integer shifts, we can without loss of generality take \(0\) to be the root.

\begin{theorem}
\label{thm:quenched-BS}
For every \(r\in\N_0\) and every finite rooted graph \(H\),
\[
\frac{1}{n}\sum_{v=1}^n\mathbf 1\{B_r(G_n,v)\cong H\}
\longrightarrow
\P\bigl(B_r(G_\infty,0)\cong H\bigr)
\qquad\text{almost surely}.
\]
Equivalently, if
\[
\nu_n=\frac1n\sum_{v=1}^n\delta_{[G_n,v]},
\qquad
\mu_\infty=\mathcal L([G_\infty,0]),
\]
then \(\nu_n\Longrightarrow\mu_\infty\) almost surely.
Here \(\Longrightarrow\) denotes Benjamini--Schramm convergence of probability measures, defined below in Definition~\ref{def:BS-convergence}.
\end{theorem}

\begin{theorem}
\label{thm:max-degree-asymptotics}
Let
\[
D_n^*=\max_{1\leq v\leq n}\deg_{G_n}(v).
\]
Then
\[
\frac{D_n^*}{\log n}\xrightarrow{\P}\frac{1}{\log(3/2)}.
\]
\end{theorem}

By Theorem~\ref{thm:degree-distribution-limit-graph}, for every integer \(d\ge2\), \(\P\bigl(\deg_{G_\infty}(0)\ge d\bigr)=(2/3)^{d-2}=e^{-(d-2)\log(3/2)}\).
Thus the coefficient \(1/\log(3/2)\) in the maximum degree asymptotics of Theorem~\ref{thm:max-degree-asymptotics} is the reciprocal of the exponential decay rate of the limiting degree distribution.

The limiting graph is almost surely one-ended and recurrent (Theorem~\ref{thm:one-endedness} and Proposition~\ref{prop:recurrence}).
These conclusions place the two-sided nearest-predecessor graph within the broader theory of recurrence for distributional limits of planar graphs and, more generally, \(H\)-minor-free graphs~\cite{GurelGurevichNachmias2013,Lee2021}.

We also obtain an exact probability generating function for the number of edges crossing the boundary of \([n]\), together with almost sure logarithmic growth and a central limit theorem.

\begin{theorem}
\label{thm:boundary-fluctuations}
For \(n\geq1\), define the \emph{edge boundary} of \([n]\) in \(G_\infty\) by
\[
\partial_E[n]
=\bigl\{\{x,y\}\in E_\infty:|\{x,y\}\cap[n]|=1\bigr\}.
\]
Then the following hold:
\begin{enumerate}
\item[(i)] For \(0<s<3/2\), the probability generating function is
\begin{equation}
\label{eq:boundary-total-pgf}
\E\left[s^{|\partial_E[n]|}\right]
=s^2\frac{\Gamma(3-2s)}{\Gamma(2s)}
\frac{\Gamma(n-1+2s)}{\Gamma(n+3-2s)}.
\end{equation}

\item[(ii)] As \(n\to\infty\),
\begin{equation}
\label{eq:boundary-as-order}
|\partial_E[n]|=\Theta(\log n)
\qquad\text{almost surely}.
\end{equation}

\item[(iii)] As \(n\to\infty\),
\begin{equation}
\label{eq:boundary-gaussian-limit}
\frac{|\partial_E[n]|-4\log n}{2\sqrt{\log n}}
\Longrightarrow\mathcal N(0,1).
\end{equation}
\end{enumerate}
\end{theorem}

\subsection{\texorpdfstring{\(K\)}{K}-Nearest-Predecessor Graph}
We next turn to the \(K\)-nearest-predecessor model.
The \(K\)-nearest-predecessor graph attaches each new vertex to its \(K\) nearest predecessors, or to all its predecessors when fewer than \(K\) are available, without prescribing how many lie on either side.

Fix an integer \(K\geq1\).
To resolve distance ties, order \(\Z\setminus\{v\}\) by
\[
v-1\prec_v v+1\prec_v v-2\prec_v v+2\prec_v\cdots,
\]
so that vertices are ordered by increasing distance from \(v\), with the left vertex first in each tie.

\begin{definition}
\label{def:finite-K-nearest-predecessor}
For \(v\in[n]\), let \(P_{n,\tau}^{(K)}(v)\) be the first \(\min\{K,\tau_n(v)-1\}\) vertices of
\[
\{u\in[n]\setminus\{v\}:\tau_n(u)<\tau_n(v)\}
\]
in the order \(\prec_v\), and define
\[
G_n^{(K)}=([n],E_n^{(K)}),
\qquad
E_n^{(K)}=\bigcup_{v\in[n]}\bigl\{\{v,u\}:u\in P_{n,\tau}^{(K)}(v)\bigr\}.
\]
\end{definition}

For every fixed \(K\), we construct the limit graph and prove the almost sure quenched convergence.

\begin{definition}
\label{def:infinite-K-nearest-predecessor}
For \(v\in\Z\), let \(P_\infty^{(K)}(v)\) be the first \(K\) vertices \(u\neq v\), in the order \(\prec_v\), satisfying \(T_u<T_v\).
These vertices exist almost surely.
Define
\[
G_\infty^{(K)}=(\Z,E_\infty^{(K)}),
\qquad
E_\infty^{(K)}=\bigcup_{v\in\Z}\bigl\{\{v,u\}:u\in P_\infty^{(K)}(v)\bigr\}.
\]
\end{definition}

\begin{theorem}
\label{thm:directionless-quenched-BS}
Let
\[
\nu_n^{(K)}=\frac1n\sum_{v=1}^n\delta_{[G_n^{(K)},v]},
\qquad
\mu_K=\mathcal L\bigl([\mathcal C_\infty^{(K)}(0),0]\bigr),
\]
where \(\mathcal C_\infty^{(K)}(0)\) is the component of \(0\) in \(G_\infty^{(K)}\).
Then, for every fixed \(K\geq1\),
\[
\nu_n^{(K)}\Longrightarrow\mu_K
\qquad\text{almost surely}.
\]
\end{theorem}

For every fixed \(K\ge1\), Lemma~\ref{lem:directionless-edge-probability} gives \(p_K(d):=\P\bigl(\{0,d\}\in E_\infty^{(K)}\bigr)\le K^2/d^2\) for every integer \(d\ge1\).
Moreover, whenever \(2d-1\ge K\), \(p_K(d)=K(K+1)/(4d^2-1)\).
Theorem~\ref{thm:directionless-one-endedness} shows that \(G_\infty^{(1)}\) is a one-ended spanning tree of \(G_\infty^{(K)}\), so \(G_\infty^{(K)}\) is connected and one-ended.
The graph \(G_\infty^{(K)}\) is also recurrent; for \(K\geq2\), it is nonplanar almost surely (Propositions~\ref{prop:directionless-recurrence} and~\ref{prop:directionless-nonplanarity}).
Thus one-endedness and recurrence persist throughout this family.

\section{The Two-Sided Nearest-Predecessor Graph}
The main goal of this section is to prove Theorem~\ref{thm:quenched-BS}, establishing quenched Benjamini--Schramm convergence for the two-sided model, which was defined in Definition~\ref{def:finite-two-sided-nearest-predecessor}.

\subsection{Degree as Record Statistics}
We first describe adjacency and degrees using record statistics.
We use standard facts about records in random permutations~\cite{Nevzorov1987}.

\begin{definition}
For \(1\leq v_1\leq v_2\leq n\), a position \(x\in[v_1,v_2]\) is a \emph{right record} if
\[
\tau_n(x)=\min_{v_1\leq y\leq x}\tau_n(y).
\]
The position \(x\) is a \emph{left record} if
\[
\tau_n(x)=\min_{x\leq y\leq v_2}\tau_n(y).
\]
Define their counts by
\[
\begin{aligned}
R_n(v_1,v_2)
&= \#\left\{x\in [v_1,v_2]:\tau_n(x)=\min_{v_1\leq y\leq x}\tau_n(y)\right\},\\
L_n(v_1,v_2)
&= \#\left\{x\in [v_1,v_2]:\tau_n(x)=\min_{x\leq y\leq v_2}\tau_n(y)\right\}.
\end{aligned}
\]
We use the convention
\[
L_n(a,b)=R_n(a,b)=0 \qquad\text{whenever }a>b.
\]
\end{definition}

\begin{lemma}[\cite{Ponomarenko2026GreedyRouting}, Edge existence criterion]
\label{lem:edge-existence-criterion}
For every \(1\leq x<y\leq n\), we have
\[
\{x,y\}\in E_n \quad\Longleftrightarrow\quad \max\{\tau_n(x),\tau_n(y)\}<\min_{x<z<y}\tau_n(z),
\]
where the minimum over an empty set is \(+\infty\).
Equivalently, \(x\) and \(y\) are adjacent if and only if they are the two earliest inserted vertices in the interval \([x,y]\).
In particular, consecutive integers \(i\) and \(i+1\) are always adjacent.
\end{lemma}

\begin{proof}
When the later of \(x\) and \(y\) is inserted, it connects to the other endpoint if and only if no vertex in \((x,y)\) has yet been inserted.
This is equivalent to both endpoints preceding every interior vertex, which is precisely the displayed condition.
The adjacent pair assertion follows because there is no integer strictly between \(i\) and \(i+1\).
\end{proof}

\begin{proposition}
\label{prop:degree-stopped-number-records}
Let \(G_n\) be the finite two-sided nearest-predecessor graph on \([n]\) with insertion time permutation \(\tau_n\).
For \(v\in[n]\), define
\[
r_n(v)=\min\bigl(\{y>v:\tau_n(y)<\tau_n(v)\}\cup\{n\}\bigr),
\]
and
\[
\ell_n(v)=\max\bigl(\{y<v:\tau_n(y)<\tau_n(v)\}\cup\{1\}\bigr).
\]
Then
\[
\deg_{G_n}(v) = L_n(\ell_n(v),v-1) + R_n(v+1,r_n(v)).
\]
\end{proposition}

\begin{proof}
We first determine the right neighbors of \(v\).
Let \(y>v\).
If \(y>r_n(v)\), then \(r_n(v)\in(v,y)\), and by definition \(\tau_n(r_n(v))<\tau_n(v)\).
Thus the edge criterion in Lemma~\ref{lem:edge-existence-criterion} fails, so \(y\) is not adjacent to \(v\).
Thus every right neighbor of \(v\) lies in \([v+1,r_n(v)]\).

Let \(v<y\le r_n(v)\).
By the definition of \(r_n(v)\), every \(z\in(v,y)\) satisfies \(\tau_n(v)<\tau_n(z)\).
Thus, if the interval is nonempty, \(\tau_n(v)<\min_{v<z<y}\tau_n(z)\); if it is empty, the inequality is vacuous.
Consequently, \(\{v,y\}\in E_n\) if and only if \(\tau_n(y)<\min_{v<z<y}\tau_n(z)\).
This is exactly the condition that \(y\) is a right record in the sequence
\[
\bigl(\tau_n(v+1),\tau_n(v+2),\ldots,\tau_n(r_n(v))\bigr).
\]
Hence, the number of right neighbors of \(v\) is \(R_n(v+1,r_n(v))\).

The symmetric argument on the left shows that the number of left neighbors is \(L_n(\ell_n(v),v-1)\).
Adding the two counts proves the formula.
\end{proof}

\subsection{The Quenched BS Limit}
\label{sec:two-sided-quenched-bs}

Theorem~\ref{thm:quenched-BS} states that the frequencies of finite rooted neighborhoods in \(G_n\) converge almost surely to their probabilities in \((G_\infty,0)\).
We begin by explaining why neighborhood frequencies give a useful notion of convergence for these sparse graphs.

Fix \(1\leq x<y\leq n\) and write \(d=y-x\).
By Lemma~\ref{lem:edge-existence-criterion}, \(x\) and \(y\) are adjacent if and only if the two endpoints receive the two smallest insertion times among the \(d+1\) vertices in \([x,y]\).
Every pair of positions in \([x,y]\) is equally likely to receive these two insertion times, so
\[
\P(\{x,y\}\in E_n)=\binom{d+1}{2}^{-1}=\frac{2}{d(d+1)}.
\]

Writing \(H_n=\sum_{k=1}^n k^{-1}\), there are \(n-d\) pairs of vertices at distance \(d\), and therefore
\[
\E\left[|E_n|\right]
=2\sum_{d=1}^{n-1}\frac{n-d}{d(d+1)}
=2\left[n\sum_{d=1}^{n-1}\left(\frac1d-\frac1{d+1}\right)-\sum_{d=1}^{n-1}\frac1{d+1}\right]
=2n-2H_n,
\]
\[
\E\left[\frac1n\sum_{v=1}^n \deg_{G_n}(v)\right]
=\frac{2\E\left[|E_n|\right]}{n}
=4-\frac{4H_n}{n}
\longrightarrow4.
\]
Thus \((G_n)_{n\geq1}\) has bounded expected average degree.

In the dense graph limit framework~\cite{LovaszSzegedy2006}, \((G_n)_{n\ge1}\)  converges in probability to the zero graphon (the limit of a sequence of edgeless graphs), while \(G_n\) remains connected with linearly many edges.
Indeed, if \(X_1,X_2\) are sampled independently and uniformly from \([n]\), then
\[
\P(\{X_1,X_2\}\in E_n) = \frac{2\E\left[|E_n|\right]}{n^2} \longrightarrow 0.
\]
Therefore, for every fixed \(k\), the subgraph of \(G_n\) induced by \(k\) uniformly sampled vertices is edgeless with probability tending to \(1\).
The zero graphon records this vanishing edge density but not the nontrivial local structure.
To study vertex degrees and neighborhood patterns, we choose a root \(\rho_n\) uniformly from \([n]\).
For each fixed radius \(r\), we examine the subgraph induced by vertices at graph distance at most \(r\) from \(\rho_n\), up to rooted isomorphism~\cite{BenjaminiSchramm2001,AldousLyons2007}.
Conditional on \(G_n\), the rooted isomorphism class has law \(\nu_n\) from Theorem~\ref{thm:quenched-BS}, with each class weighted by the proportion of vertices that produce it when chosen as the root.

The next three definitions introduce a topology on rooted graphs, weak convergence of their laws, and almost sure weak convergence of these empirical measures.

\begin{definition}
A rooted graph is a graph \(G\) with a distinguished vertex \(o\), and \(B_r(G,o)\) is the rooted subgraph induced by the vertices within graph distance \(r\) of \(o\).
Let \(\mG_*\) be the \emph{space of rooted isomorphism classes} of connected, locally finite rooted graphs. 
Define
\[
q=\sup\{r\in\N_0:B_r(G,o)\cong B_r(G',o')\},
\qquad
d_{\mathrm{loc}}\bigl((G,o),(G',o')\bigr)=\frac{1}{1+q}.
\]
\end{definition}

The open balls of \(d_{\mathrm{loc}}\) generate the \emph{local topology} on \(\mG_*\)~\cite{GurelGurevichNachmias2013}.
Thus \((G_m,o_m)\to(G,o)\) in the local topology if and only if, for every fixed \(r\), the balls \(B_r(G_m,o_m)\) and \(B_r(G,o)\) are isomorphic for all sufficiently large \(m\).

\begin{definition}
\label{def:BS-convergence}
The weak topology of probability measures on \((\mG_*,d_{\mathrm{loc}})\) is called the \emph{local weak topology}.
Convergence in this topology is called Benjamini--Schramm (BS) convergence, or \emph{local weak convergence}.

\end{definition}

\begin{definition}
Let \(G_n\) be a random graph on \([n]\).
Here \([G,v]\) denotes the rooted isomorphism class of \((G,v)\).
We say that \(G_n\) converges to a random rooted graph \((G,o)\) in the \emph{quenched Benjamini--Schramm sense} if
\[
\nu_n:=\frac1n\sum_{v=1}^n\delta_{[G_n,v]}
\Longrightarrow \mathcal L([G,o])
\qquad\text{almost surely},
\]
where \(\Longrightarrow\) denotes convergence in the local weak topology (Definition~\ref{def:BS-convergence}).
\end{definition}

Recall the infinite graph from Definition~\ref{def:infinite-two-sided-nearest-predecessor}.

\begin{lemma}
\label{lem:local-finiteness-of-the-infinite-graph}
The graph \(G_\infty\) is locally finite almost surely.
\end{lemma}

\begin{proof}
Set
\[
N_+=\inf\{j\ge1:T_j<T_0\},
\qquad
N_-=\inf\{j\ge1:T_{-j}<T_0\}.
\]
For every \(m\ge1\),
\[
\P(N_+>m)=\P(N_->m)=\int_0^1(1-t)^m\,dt=\frac1{m+1},
\]
so \(N_+,N_-<\infty\) almost surely.
If \(j>N_+\), then
\[
\min_{0<k<j}T_k\le T_{N_+}<T_0\le\max\{T_0,T_j\},
\]
so \(j\) is not a neighbor of \(0\).
The left side is analogous, so all neighbors of \(0\) lie in \(\{-N_-,\ldots,N_+\}\).

Thus \(\deg_{G_\infty}(0)\le N_++N_-<\infty\) almost surely.
By translation invariance, \(\P(\deg_{G_\infty}(v)=\infty)=0\) for every \(v\in\Z\), and hence
\[
\P\left(\exists v\in\Z:\deg_{G_\infty}(v)=\infty\right)
\leq\sum_{v\in\Z}\P\left(\deg_{G_\infty}(v)=\infty\right)
=0.
\]
\end{proof}

The following lemma links the two equivalent assertions of Theorem~\ref{thm:quenched-BS}.

\begin{lemma}
\label{lem:cylinder-events-determine-BS-convergence}
For \(r\in\N_0\) and a finite rooted graph \(H\), define \emph{cylinder events}
\[
C_{r,H} = \{(G,o)\in\mG_*:B_r(G,o)\cong H\}.
\]
Then the family of events \(C_{r,H}\) is countable and convergence-determining.
That is, if \((\nu_n)\) and \(\mu\) are probability measures on \(\mG_*\) such that, for every \(r\in\N_0\) and every finite rooted graph \(H\),
\begin{equation}
\label{eq:cylinder-convergence}
\nu_n(C_{r,H})\longrightarrow \mu(C_{r,H}),
\end{equation}
then Benjamini--Schramm convergence holds:
\[
\nu_n\Longrightarrow \mu.
\]
\end{lemma}

\noindent The proof is given in Appendix~\ref{app:proof-cylinder-events}.

To prove Theorem~\ref{thm:quenched-BS}, we first couple the finite and infinite graphs.

\begin{lemma}
\label{lem:projective-consistency}
For each \(n\geq1\), define
\begin{equation}
\label{eq:finite-rank-permutation}
\tau_n(i)=1+\sum_{j=1}^n \mathbf 1\{T_j<T_i\},
\qquad 1\leq i\leq n.
\end{equation}
Then \(\tau_n\) is a uniformly random permutation of \([n]\), and
\[
G_n=G_\infty[[n]]
\qquad\text{almost surely for every }n\geq1.
\]
\end{lemma}

\begin{proof}
Since \(T_1,\ldots,T_n\) are i.i.d.\ continuous random variables, their relative ranks form a uniformly random permutation of \([n]\).

Fix \(n\ge1\) and \(1\leq i<j\leq n\).
By Lemma~\ref{lem:edge-existence-criterion}, the edge \(\{i,j\}\) belongs to \(G_n\) if and only if \(\max\{\tau_n(i),\tau_n(j)\} < \min_{i<k<j}\tau_n(k)\).
By~\eqref{eq:finite-rank-permutation}, this is equivalent to \(\max\{T_i,T_j\}<\min_{i<k<j}T_k\).
This is exactly the condition that \(\{i,j\}\in E_\infty\).
Hence the edge sets agree on all pairs inside \([n]\).
Since there are only countably many choices of \(n,i,j\), the equality holds almost surely for every \(n\).
\end{proof}

We call this compatibility under restriction \emph{projective consistency}: under the coupling, restricting \(G_\infty\), or any \(G_m\) with \(m\ge n\), to \([n]\) gives \(G_n\) almost surely.

Moving the root \(v\) to \(0\) by translation is equivalent, up to rooted isomorphism, to shifting the insertion-time sequence by \(v\).
Fact~\ref{fact:shift-ergodic} supplies the ergodicity needed to apply Birkhoff's theorem to neighborhood frequencies.
\begin{fact}[{\cite{Shields1973BernoulliShifts}}]
\label{fact:shift-ergodic}
Let
\[
\Omega=[0,1]^{\Z}, \qquad \mathcal F=\mathcal B([0,1])^{\otimes\Z}, \qquad \P=\lambda^{\otimes\Z},
\]
where \(\lambda\) is Lebesgue measure on \([0,1]\).
Define the shift
\[
\theta:\Omega\to\Omega, \qquad (\theta\omega)_i=\omega_{i+1}.
\]
Then \(\theta\) is measure preserving and ergodic.
\end{fact}

\begin{proof}[Proof of Theorem~\ref{thm:quenched-BS}]
We first prove the neighborhood-frequency limits in the first assertion, and then deduce the second assertion using Lemma~\ref{lem:cylinder-events-determine-BS-convergence}.

Fix \(r\in\N_0\) and a finite rooted graph \(H\).
For insertion time \(\omega=(T_i)_{i\in\Z}\), define the shift and the local statistic together by
\[
(\theta\omega)_i=\omega_{i+1},\qquad
X_v(\omega):=\mathbf 1\{B_r(G_\infty(\omega),v)\cong H\}
=X_0(\theta^v\omega).
\]
Fact~\ref{fact:shift-ergodic} and Birkhoff's ergodic theorem yield
\[
\frac1n\sum_{v=1}^n X_v = \frac1n\sum_{v=1}^n X_0\circ\theta^v \longrightarrow \E\left[X_0\right] = \P\left(B_r(G_\infty,0)\cong H\right) \qquad\text{a.s.}
\]
Set \(Y_{n,v}=\mathbf 1\{B_r(G_n,v)\cong H\}\). Then
\[
\left| \frac1n\sum_{v=1}^n Y_{n,v} - \frac1n\sum_{v=1}^n X_v \right| \leq \frac1n \sum_{v=1}^n \mathbf 1\{B_r(G_n,v)\neq B_r(G_\infty,v)\}.
\]
For \(v\in\Z\), define
\[
M_r(v)=\max\{|x-v|:x\in B_r(G_\infty,v)\},
\qquad d_n(v)=\min\{v-1,n-v\}.
\]
By Lemma~\ref{lem:local-finiteness-of-the-infinite-graph}, \(M_r(v)<\infty\) almost surely.
If \(M_r(v)\le d_n(v)\), then \(B_r(G_\infty,v)\subseteq[n]\), and every shortest path from \(v\) to a vertex of this ball stays inside \([n]\).
Since \(G_n=G_\infty[[n]]\) by Lemma~\ref{lem:projective-consistency}, these paths and all induced edges are preserved in \(G_n\), while every path in \(G_n\) is also a path in \(G_\infty\).
Therefore \(B_r(G_n,v)=B_r(G_\infty,v)\). Hence
\[
\frac1n \sum_{v=1}^n \mathbf 1\{B_r(G_n,v)\neq B_r(G_\infty,v)\} \leq \frac1n \sum_{v=1}^n \mathbf 1\{M_r(v)>d_n(v)\}.
\]
Fix \(h\in\N\).
The vertices satisfying \(d_n(v)<h\) lie in the two boundary intervals
\[
\{1,\ldots,h\}
\qquad\text{and}\qquad
\{n-h+1,\ldots,n\},
\]
whose total size is at most \(2h\).
Hence
\[
\frac1n \sum_{v=1}^n \mathbf 1\{M_r(v)>d_n(v)\} \leq \frac{2h}{n} + \frac1n \sum_{v=h+1}^{n-h} \mathbf 1\{M_r(v)>h\}.
\]
Since \(M_r(v)=M_r(0)\circ\theta^v\), another application of Birkhoff's theorem gives, for every \(h\in\N\),
\[
\frac1n \sum_{v=1}^n \mathbf 1\{M_r(v)>h\} \longrightarrow \P(M_r(0)>h) \qquad\text{a.s.}
\]
Therefore, for every \(h\in\N\),
\[
\limsup_{n\to\infty} \left| \frac1n\sum_{v=1}^n Y_{n,v} - \frac1n\sum_{v=1}^n X_v \right| \leq \P(M_r(0)>h) \qquad\text{a.s.}
\]
Letting \(h\to\infty\), using \(M_r(0)<\infty\) almost surely, yields
\[
\frac1n \sum_{v=1}^n \mathbf 1\{B_r(G_n,v)\cong H\} - \frac1n \sum_{v=1}^n \mathbf 1\{B_r(G_\infty,v)\cong H\} \longrightarrow 0 \qquad\text{a.s.}
\]
Combining this with the Birkhoff limit for the second term yields, for the fixed pair \((r,H)\),
\[
\frac1n \sum_{v=1}^n \mathbf 1\{B_r(G_n,v)\cong H\} \longrightarrow \P\left(B_r(G_\infty,0)\cong H\right) \qquad\text{a.s.}
\]
Since there are only countably many pairs \((r,H)\), the intersection of the events on which the corresponding convergences hold has probability one.
Therefore, almost surely,
\begin{equation}
\label{eq:rooted-ball-measure-identities}
\nu_n(C_{r,H})
=\frac1n\sum_{v=1}^n\mathbf 1\{B_r(G_n,v)\cong H\}
\longrightarrow \P\bigl(B_r(G_\infty,0)\cong H\bigr)
=\mu_\infty(C_{r,H})
\end{equation}
for every \((r,H)\).
Applying Lemma~\ref{lem:cylinder-events-determine-BS-convergence} gives \(\nu_n\Longrightarrow\mu_\infty\) almost surely, proving the second assertion.
The reverse implication follows because each cylinder event \(C_{r,H}\) is both open and closed in the local topology.
\end{proof}

\begin{remark}
Theorem~\ref{thm:quenched-BS} uniquely determines the law of (the isomorphism class of) the rooted graph \((G_\infty,0)\).
This is because the space of measures on \(\mG_*\) is Hausdorff under the local weak topology.
\end{remark}

\begin{remark}
\label{rem:quenched-implies-annealed}
Let \(\rho_n\) be uniform on \([n]\), independently of the insertion times.
For every \(r\in\N_0\) and every finite rooted graph \(H\), conditioning on \(G_n\) and then applying Theorem~\ref{thm:quenched-BS} and bounded convergence gives
\[
\P\bigl(B_r(G_n,\rho_n)\cong H\bigr)
=
\E\left[\frac1n\sum_{v=1}^n
\mathbf 1\{B_r(G_n,v)\cong H\}\right]
\longrightarrow
\P\bigl(B_r(G_\infty,0)\cong H\bigr).
\]
This convergence of the uniformly rooted graphs after averaging over the insertion times is called \emph{annealed Benjamini--Schramm convergence}, and we write
\[
(G_n,\rho_n)\xrightarrow{\mathrm{BS}}(G_\infty,0).
\]
\end{remark}

\section{Structural Properties of the Two-Sided Nearest-Predecessor Graph}
The rooted graph \((G_\infty,0)\) describes the limiting local neighborhood of a uniformly chosen vertex of \(G_n\).
We use this limit graph as the main object for computing and organizing structural information, and then transfer the relevant finite graph consequences back through projective consistency and local convergence.

\subsection{The limiting degree distribution}
\begin{theorem}
\label{thm:degree-distribution-limit-graph}
Let \(D=\deg_{G_\infty}(0)\).
Then, for every integer \(d\ge2\),
\[
\P(D=d) = \frac13\left(\frac23\right)^{d-2}.
\]
In particular, \(D\) has an exponential tail.\footnote{Luque, Lacasa, Ballesteros, and Luque proved the same law for horizontal visibility graphs by enumerating visibility configurations and evaluating nested integrals~\cite{LuqueLacasaBallesterosLuque2010}.}
\end{theorem}

\begin{proof}
Let \(D_+\) and \(D_-\) be the numbers of neighbors of \(0\) on the two sides, and set \(A:=-\log T_0\sim\operatorname{Exp}(1)\).
Conditional on \(A=a\), the transformed insertion times \(Y_j=-\log T_j\) on either half line are i.i.d.\ \(\operatorname{Exp}(1)\), and~\eqref{inf-rule} gives, on the right,
\[
\{0,j\}\in E_\infty \quad\Longleftrightarrow\quad
\min\{a,Y_j\}>\max_{1\le k<j}Y_k.
\]
Let \(R_1<R_2<\cdots\) denote the successive record values of the sequence \((Y_j)\).
Conditional on the record history up to \(R_k=x\), the next record is the first subsequent observation exceeding \(x\).
By the memoryless property of the exponential distribution,
\[
R_{k+1}-R_k\sim\operatorname{Exp}(1),
\]
independently of the preceding record values.
The initial record \(R_1\) also has distribution \(\operatorname{Exp}(1)\).
Hence the record values are partial sums of independent rate-one exponential variables and therefore are the arrival times of a rate-one Poisson process.
By the same argument as in the proof of Proposition~\ref{prop:degree-stopped-number-records}, applied to the infinite edge rule~\eqref{inf-rule}, the right neighbors of \(0\) correspond to the record values not exceeding \(a\), together with the first record value exceeding \(a\).
The left side is analogous and, conditional on \(A=a\), independent of the right side.
Consequently,
\[
D_+-1,D_--1\mid(A=a)\stackrel{\mathrm{ind}}{\sim}\operatorname{Poisson}(a),
\qquad
D-2\mid(A=a)\sim\operatorname{Poisson}(2a).
\]
Integrating over \(A\), for every integer \(d\ge2\),
\[
\P(D=d) = \int_0^\infty e^{-2a}\frac{(2a)^{d-2}}{(d-2)!} e^{-a}\,da
= \frac{2^{d-2}}{(d-2)!} \cdot \frac{(d-2)!}{3^{d-1}} = \frac13\left(\frac23\right)^{d-2}.
\]
\end{proof}

\subsection{The Maximum Degree Asymptotics}
We prove Theorem~\ref{thm:max-degree-asymptotics} by obtaining an upper bound from the limiting degree law and a matching lower bound from record counts in independent blocks of vertices.

\begin{lemma}
\label{lem:record-generating-function}
Let \(m\ge 0\), and let \(R_m\) be the number of right records in a uniform random permutation on \(m\) elements.
Then
\[
C_m(z):=\E\left[z^{R_m}\right] =
\begin{cases}
1, & m=0,\\[3pt]
\dfrac{z(z+1)\cdots(z+m-1)}{m!}, & m\ge 1.
\end{cases}
\]
For \(z>0\) and \(m\ge1\), this can also be written as
\[
C_m(z)=\frac{\Gamma(m+z)}{\Gamma(z)\Gamma(m+1)}.
\]
The same generating function holds for the number of left records.
\end{lemma}

\noindent The proof is given in Appendix~\ref{app:proof-record-generating-function}.

\begin{fact}[{\cite{TricomiErdelyi1951}}]
\label{fact:gamma-ratio}
Let \(J\subset\R\) be compact.
Uniformly for \(a,b\in J\),
\[
\frac{\Gamma(n+a)}{\Gamma(n+b)}
=n^{a-b}\left(1+O_J\left(n^{-1}\right)\right)
\qquad\text{as }n\to\infty.
\]
\end{fact}

\begin{lemma}
\label{lem:two-sided-number-records-estimate}
For \(m\ge 1\), let \(L_m\) and \(R_m\) be independent random variables distributed as the numbers of left records and right records, respectively, in a uniform random permutation on \(m\) elements.
Both have probability generating function (PGF) \(C_m(z)\).
Define
\[
S_m:=L_m+R_m.
\]
For \(u>2\), define
\[
I(u):=u\log\frac u2-u+2.
\]
Then, as \(m\to\infty\),
\[
\P(S_m\ge u\log m) = m^{-I(u)+o(1)}.
\]
\end{lemma}

\noindent The proof, given in Appendix~\ref{app:proof-two-sided-number-records}, uses \emph{exponential tilting}, a change of probability measure that gives exponentially greater weight to larger values of \(S_m\), so that the upper tail event has probability tending to \(1\) under the tilted measure.
The tilted measure is defined explicitly in~\eqref{measure-Q}.

\begin{proof}[Proof of Theorem~\ref{thm:max-degree-asymptotics}]
We want to show that, for every \(\varepsilon>0\),
\[
\P\left(
(1-\varepsilon)\frac{\log n}{\log(3/2)}
\leq D_n^*
\leq
(1+\varepsilon)\frac{\log n}{\log(3/2)}
\right)
\longrightarrow 1.
\]
We first prove the upper bound.
By Lemma~\ref{lem:projective-consistency}, every edge incident to \(v\) in \(G_n\) is also incident to \(v\) in \(G_\infty\).
\(\deg_{G_n}(v)\le D_v:=\deg_{G_\infty}(v)\).
So, by Theorem~\ref{thm:degree-distribution-limit-graph},
\[
\P(D_v\ge d)=(2/3)^{d-2}, \qquad d \ge 2.
\]
Hence the union bound gives
\[
\P(D_n^*\ge d)\le \sum_{v=1}^n\P(D_v\ge d)
=n\left(\frac23\right)^{d-2}.
\]
Taking \(d_n=\lceil(1+\varepsilon)\log n/\log(3/2)\rceil\),
\[
\P(D_n^*\ge d_n) \le n\left(\frac23\right)^{d_n-2} \le \left(\frac32\right)^2 n^{-\varepsilon} \longrightarrow 0.
\]
We now prove the lower bound.
We may assume \(0<\varepsilon<1\), since for \(\varepsilon\ge 1\) the lower bound is trivial.
Partition the first \(q w_n\) vertices into blocks of width \(w_n\), where
\[
\alpha:=\frac{1-\varepsilon}{3\log(3/2)},
\qquad
w_n := \left\lfloor n^\alpha \right\rfloor,
\qquad q:=\left\lfloor\frac n{w_n}\right\rfloor.
\]

Let \(Y_b\) be the within-block degree of the earliest vertex in block \(b\).
The \(Y_b\)'s are independent because they depend only on the relative orders in disjoint blocks, and their edges remain present in \(G_n\); hence
\[
D_n^*\ge \max_{1\le b\le q}Y_b.
\]
We estimate
\[
p_n := \P\left( Y_b \ge (1-\varepsilon)\frac{\log n}{\log(3/2)} \right).
\]
Let \(v_b^\star\) be the earliest-inserted vertex inside the \(b\)-th block, and let \(J_b\in\{1,\ldots,w_n\}\) be its position relative to the left endpoint of that block.
Then \(J_b\) is uniform, and with
\[
\mathcal M_b=\left\{\left\lceil\frac{w_n}3\right\rceil\le J_b\le\left\lfloor\frac{2w_n}{3}\right\rfloor\right\},
\qquad m=\left\lfloor\frac{w_n}6\right\rfloor,
\]
we have \(\P(\mathcal M_b)=1/3+O(1/w_n)\).
On \(\mathcal M_b\), let \(I_b^-\) and \(I_b^+\) be the \(m\) vertices immediately to the left and right of \(v_b^\star\); both lie in the block for large \(w_n\).
Because \(v_b^\star\) is the earliest block vertex and these intervals are adjacent to it, Lemma~\ref{lem:edge-existence-criterion} and independence of their relative orders give
\[
Y_b\ge L_m^{(b)}+R_m^{(b)}\quad\text{on }\mathcal M_b,
\qquad
\bigl(L_m^{(b)}+R_m^{(b)}\mid\mathcal M_b\bigr)\stackrel d=L_m+R_m=S_m,
\]
where \(L_m\) and \(R_m\) are independent record counts on \(m\) elements.

Since \(m=\lfloor w_n/6\rfloor\), we have \(\log m=\alpha\log n+O(1)\), and hence
\[
a_m := \frac{(1-\varepsilon)\log n} {\log(3/2)\log m} =3+O\left(\frac1{\log m}\right).
\]
For every fixed \(\delta\in(0,1)\), all sufficiently large \(m\) satisfy \(3-\delta\leq a_m\leq3+\delta\). Hence
\[
\P(S_m\geq(3+\delta)\log m) \leq \P(S_m\geq a_m\log m) \leq \P(S_m\geq(3-\delta)\log m).
\]
Applying Lemma~\ref{lem:two-sided-number-records-estimate} and then letting \(\delta\downarrow0\) gives
\[
\P(S_m\geq a_m\log m) = m^{-I(3)+o(1)}.
\]
Since \(I(3)=3\log\frac32-1\), we have
\[
p_n \ge \left(\frac13+o(1)\right) \P\left(S_m \ge a_m\log m \right) = \left(\frac13+o(1)\right) m^{-I(3)+o(1)} = w_n^{-I(3)+o(1)}.
\]
Therefore,
\[
\begin{aligned}
q p_n &\ge \left(\frac n{w_n}+O(1)\right)w_n^{-I(3)+o(1)}
=nw_n^{-1-I(3)+o(1)}\\[4pt]
&=n^{1-\alpha(1+I(3))+o(1)}
=n^{1-\alpha 3\log(3/2)+o(1)}
=n^{\varepsilon+o(1)}.
\end{aligned}
\]
Thus \(q p_n\to\infty\).
By the independence of \(Y_1,\ldots,Y_q\),
\[
\begin{aligned}
\P\left( D_n^* < (1-\varepsilon)\frac{\log n}{\log(3/2)} \right) &\le \P\left( \max_{1\le b\le q}Y_b < (1-\varepsilon)\frac{\log n}{\log(3/2)} \right)\\[5pt]
&= (1-p_n)^q \le \exp(-q p_n)\\[2pt]
&\longrightarrow 0.
\end{aligned}
\]
\end{proof}

\subsection{One-Endedness and Recurrence}

\begin{definition}
Let \(G\) be an infinite, connected, locally finite graph.
An \emph{infinite simple ray} is a one-way infinite simple path in \(G\).
An \emph{end} of \(G\) is an equivalence class of infinite simple rays, where two rays are equivalent if, after deleting any finite set of vertices, sufficiently long tails of the two rays lie in the same connected component.

The graph \(G\) is \emph{one-ended} if it has exactly one end; equivalently, \(G\setminus S\) has one infinite connected component for every finite \(S\subset V(G)\).
\end{definition}

\begin{theorem}
\label{thm:one-endedness}
The graph \(G_\infty\) is almost surely one-ended.
\end{theorem}

\begin{proof}
Fix a finite interval \([a,b]\) and let
\[
\kappa=\arg\min_{a\le k\le b}T_k,
\]
which is unique almost surely.
Since \(T_\kappa>0\), insertion times below \(T_\kappa\) appear infinitely often on each exterior half line, so almost surely
\[
x_-:=\max\{\ell<a:T_\ell<T_\kappa\},\qquad
x_+:=\min\{r>b:T_r<T_\kappa\}
\]
exist.
By their choice,
\[
\max\{T_{x_-},T_{x_+}\}<T_\kappa
=\min_{x_-<k<x_+}T_k,
\]
and the edge criterion~\eqref{inf-rule} gives \(\{x_-,x_+\}\in E_\infty\).
Since the integer intervals are countable, this bypass construction holds for every \([a,b]\) on a single event of probability one.

On this event, let \(S\subset\Z\) be finite, choose \([a,b]\supseteq S\), and let \(x_-,x_+\) be the corresponding bypass vertices.
The nearest-neighbor paths on \((-\infty,a-1]\) and \([b+1,\infty)\), together with this bypass edge, survive in \(G_\infty\setminus S\).
They form one infinite component containing every vertex outside \([a,b]\); any other component is contained in \([a,b]\) and is therefore finite.
Thus \(G_\infty\setminus S\) has exactly one infinite component, proving one-endedness.
\end{proof}

\begin{definition}
Let \(G\) be an infinite, connected, locally finite graph, and let \((X_m)_{m\ge0}\) be a simple random walk on \(G\).
We call \(G\) \emph{recurrent} if, for one (hence every) \(o\in V(G)\),
\[
\tau_o^+=\inf\{m\geq1:X_m=o\}, \qquad \P_o(\tau_o^+<\infty)=1.
\]
\end{definition}

\begin{theorem}[{Gurel-Gurevich and Nachmias~\cite[Theorem~1.1]{GurelGurevichNachmias2013}}]
\label{thm:planar-limit-recurrence}
Let \((G_n,\rho_n)\) be finite planar graphs rooted at a uniformly chosen vertex, and suppose \((G_n,\rho_n)\xrightarrow{\mathrm{BS}}(G,o)\).
If \((G,o)\) is almost surely infinite and the degree of \(o\) in \(G\) has an exponential tail, then \(G\) is almost surely recurrent.
\end{theorem}

\begin{proposition}
\label{prop:recurrence}
The graph \(G_\infty\) is almost surely recurrent.
\end{proposition}

\begin{proof}
Remark~\ref{rem:quenched-implies-annealed} gives annealed BS convergence to \((G_\infty,0)\).
Each \(G_n\) is outerplanar, hence planar, because placing its vertices in their natural order on a circle makes crossing chords \(a<c<b<d\) impossible: Lemma~\ref{lem:edge-existence-criterion} would force both \(\tau_n(c)>\tau_n(b)\) and \(\tau_n(b)>\tau_n(c)\).
Theorem~\ref{thm:degree-distribution-limit-graph} gives an exponential tail for the root degree.
Thus, the claim follows from Theorem~\ref{thm:planar-limit-recurrence}.
\end{proof}

\subsection{Edge Boundaries}
Our goal is to prove Theorem~\ref{thm:boundary-fluctuations} by expressing the number of boundary edges in terms of record counts.

Here we introduce notation to be used in the following lemma and the proof of Theorem~\ref{thm:boundary-fluctuations}.
Fix \(n\ge1\), let \(B_n^+\) and \(B_n^-\) count the right-crossing and left-crossing boundary edges of \([n]\), respectively, and set
\[
X_n:=|\partial_E[n]|=B_n^-+B_n^+.
\]
List the left records in \([n]\) and the right records beyond it as
\[
\begin{aligned}
\ell_1&=n>\cdots>\ell_{L_n}=\kappa_n,
&\ell_{p+1}&=\max\{i<\ell_p:T_i<T_{\ell_p}\},\\
r_1&=n+1<r_2<\cdots,
&r_{q+1}&=\min\{j>r_q:T_j<T_{r_q}\}.
\end{aligned}
\]
Write \(a_p=T_{\ell_p}\), \(u_q=T_{r_q}\), and \(m_n=a_{L_n}=\min_{i\in[n]}T_i\).
Since all insertion times lie in \((0,1)\) almost surely, set \(a_0=u_0=1\) as a common endpoint above both record sequences.
Let \(Q=R_n^+(m_n):=\#\{q\ge1:u_q>m_n\}\).

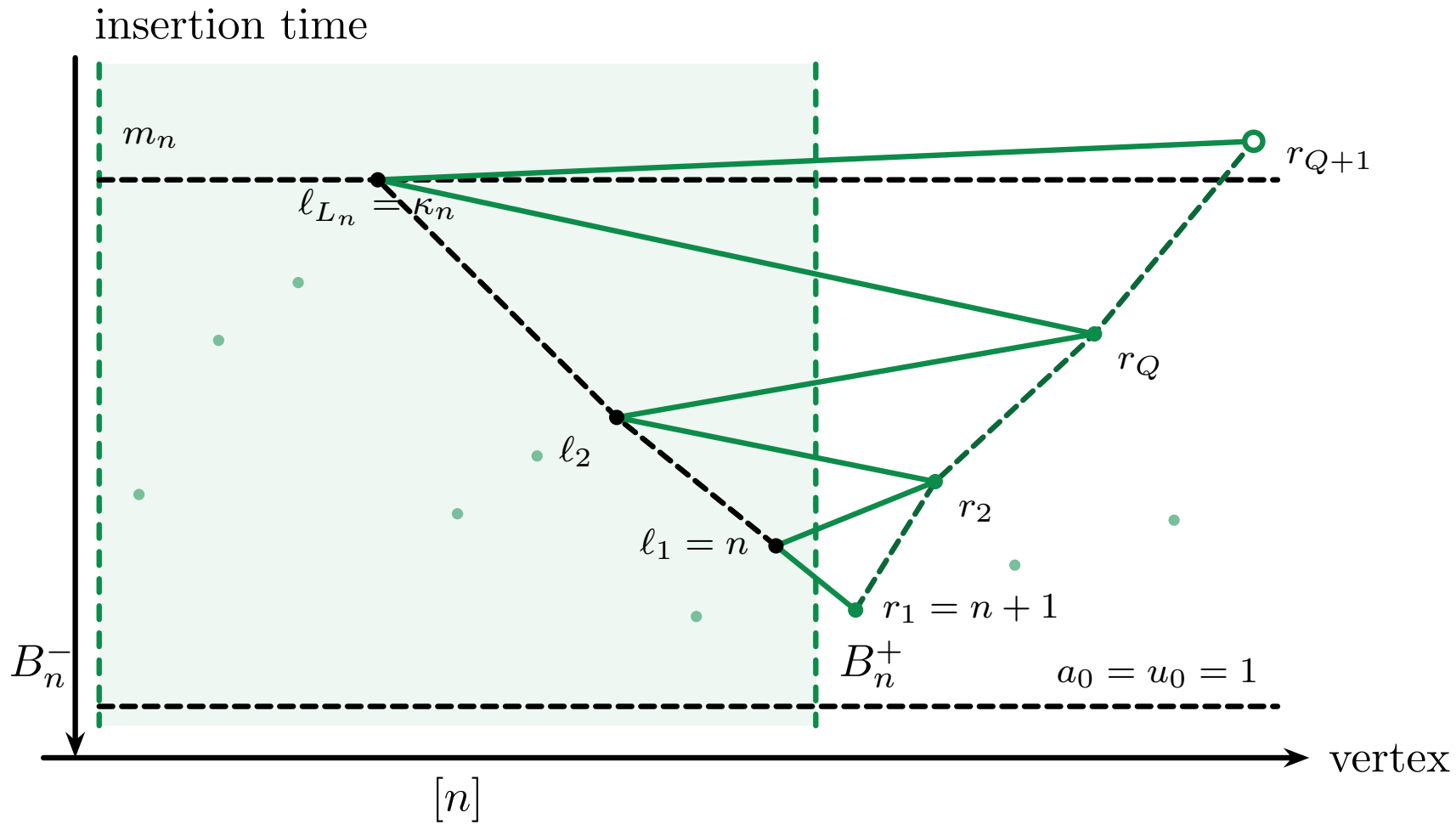
\begin{figure}[ht]
\centering
\begin{tikzpicture}[x=0.62cm,y=-5.0cm,line width=1.2pt,line cap=round,line join=round]
\fill[babyblue!7] (0.5,0) rectangle (9.5,1.03);
\draw[black,-{Stealth[length=2mm]}] (-0.2,1.08) -- (15.7,1.08)
    node[right,font=\small] {vertex};
\draw[black,-{Stealth[length=2mm]}] (0.2,-0.01) -- (0.2,1.08);
\node[anchor=south west,font=\small] at (0.2,-0.01) {insertion time};
\draw[babyblue,dashed] (0.5,0) -- (0.5,1.03);
\draw[babyblue,dashed] (9.5,0) -- (9.5,1.03);
\node[font=\small] at (5,1.145) {\([n]\)};
\node[black,font=\small,anchor=south east] at (0.45,1) {\(B_n^-\)};
\node[black,font=\small,anchor=south west] at (9.55,1) {\(B_n^+\)};

\draw[black,densely dashed] (0.5,1) -- (15.3,1);
\node[font=\scriptsize,anchor=south east]
    at (15.3,1) {\(a_0=u_0=1\)};
\draw[black,densely dashed] (0.5,0.18) -- (15.3,0.18);
\node[font=\scriptsize,anchor=south west]
    at (0.55,0.16) {\(m_n\)};

\foreach \x/\y in {1/0.67,2/0.43,3/0.34,5/0.70,6/0.61,8/0.86,
                     12/0.78,14/0.71} {
    \fill[babyblue!55] (\x,\y) circle (1.25pt);
}

\draw[black,densely dashed] (9,0.75) -- (7,0.55) -- (4,0.18);
\draw[babyblue!75!black,densely dashed]
    (10,0.85) -- (11,0.65) -- (13,0.42) -- (15,0.12);

\foreach \xa/\ya/\xb/\yb in {
    9/0.75/10/0.85,
    9/0.75/11/0.65,
    7/0.55/11/0.65,
    7/0.55/13/0.42,
    4/0.18/13/0.42,
    4/0.18/15/0.12
} {
    \draw[babyblue,line width=1.2pt] (\xa,\ya) -- (\xb,\yb);
}

\foreach \x/\y in {9/0.75,7/0.55,4/0.18} {
    \fill[black] (\x,\y) circle (1.7pt);
}
\foreach \x/\y in {10/0.85,11/0.65,13/0.42} {
    \fill[babyblue] (\x,\y) circle (1.7pt);
}
\fill[white] (15,0.12) circle (2pt);
\draw[babyblue,line width=1.2pt] (15,0.12) circle (2pt);

\node[font=\scriptsize,anchor=south east]
    at (8.9,0.8) {\(\ell_1=n\)};
\node[font=\scriptsize,anchor=north east]
    at (6.95,0.55) {\(\ell_2\)};
\node[font=\scriptsize,anchor=north]
    at (4,0.165) {\(\ell_{L_n}=\kappa_n\)};
\node[font=\scriptsize,anchor=south west]
    at (10.1,0.9) {\(r_1=n+1\)};
\node[font=\scriptsize,anchor=north west]
    at (11.05,0.65) {\(r_2\)};
\node[font=\scriptsize,anchor=north west]
    at (13.05,0.42) {\(r_Q\)};
\node[font=\scriptsize,anchor=north east]
    at (16.75,0.1) {\(r_{Q+1}\)};
\end{tikzpicture}
\caption{The solid blue segments are the right-crossing edges counted by \(B_n^+\); the hollow point marks the first right record \(r_{Q+1}\) with \(u_{Q+1}<m_n\).}
\label{fig:boundary-records}
\end{figure}

\begin{lemma}
\label{lem:boundary-common-refinement}
For \(1\le p\le L_n\) and \(1\le q\le Q+1\), set
\[
A_p=[a_p,a_{p-1}],
\qquad
U_q=[\max\{m_n,u_q\},u_{q-1}].
\]
Then \((A_p)\) and \((U_q)\) are partitions of \([m_n,1]\), and the right-crossing edges are in bijection with the nondegenerate cells in their \emph{common refinement}.
\end{lemma}

\begin{figure}[ht]
\centering
\begin{tikzpicture}[x=1cm,y=-5.8cm,line cap=round,line join=round]
\node[font=\small] at (1,0.10) {\((A_p)\)};
\node[font=\small] at (3.5,0.10) {\((U_q)\)};
\node[font=\small] at (7.7,0.10) {common refinement};

\foreach \y in {0.18,0.42,0.55,0.65,0.75,0.85,1} {
    \draw[babyblue!35,densely dashed,line width=0.8pt]
        (0.65,\y) -- (6.35,\y);
}

\draw[black,line width=1.6pt] (1,0.18) -- (1,0.55);
\draw[black,line width=1.6pt] (1,0.55) -- (1,0.75);
\draw[black,line width=1.6pt] (1,0.75) -- (1,1);
\foreach \y in {0.18,0.55,0.75,1} {
    \draw[black,line width=1.2pt] (0.88,\y) -- (1.12,\y);
}
\node[font=\scriptsize,anchor=east] at (0.85,0.365) {\(A_3\)};
\node[font=\scriptsize,anchor=east] at (0.85,0.65) {\(A_2\)};
\node[font=\scriptsize,anchor=east] at (0.85,0.875) {\(A_1\)};

\draw[babyblue,line width=1.6pt] (3.5,0.18) -- (3.5,0.42);
\draw[babyblue,line width=1.6pt] (3.5,0.42) -- (3.5,0.65);
\draw[babyblue,line width=1.6pt] (3.5,0.65) -- (3.5,0.85);
\draw[babyblue,line width=1.6pt] (3.5,0.85) -- (3.5,1);
\foreach \y in {0.18,0.42,0.65,0.85,1} {
    \draw[babyblue,line width=1.2pt] (3.38,\y) -- (3.62,\y);
}
\node[font=\scriptsize,anchor=west] at (3.65,0.30) {\(U_4\)};
\node[font=\scriptsize,anchor=west] at (3.65,0.535) {\(U_3\)};
\node[font=\scriptsize,anchor=west] at (3.65,0.75) {\(U_2\)};
\node[font=\scriptsize,anchor=west] at (3.65,0.925) {\(U_1\)};

\foreach \ya/\yb/\ym/\p/\q in {
    0.18/0.42/0.30/3/4,
    0.42/0.55/0.485/3/3,
    0.55/0.65/0.60/2/3,
    0.65/0.75/0.70/2/2,
    0.75/0.85/0.80/1/2,
    0.85/1/0.925/1/1
} {
    \draw[babyblue,line width=2.2pt] (6.4,\ya) -- (6.4,\yb);
    \node[font=\scriptsize,anchor=west] at (6.58,\ym)
        {\(A_{\p}\cap U_{\q}\longleftrightarrow\{\ell_{\p},r_{\q}\}\)};
}
\foreach \y in {0.18,0.42,0.55,0.65,0.75,0.85,1} {
    \draw[black,line width=1.2pt] (6.28,\y) -- (6.52,\y);
}

\node[font=\scriptsize,anchor=east] at (6.18,0.18) {\(m_n\)};
\node[font=\scriptsize,anchor=east] at (6.18,0.42) {\(u_3\)};
\node[font=\scriptsize,anchor=east] at (6.18,0.55) {\(a_2\)};
\node[font=\scriptsize,anchor=east] at (6.18,0.65) {\(u_2\)};
\node[font=\scriptsize,anchor=east] at (6.18,0.75) {\(a_1\)};
\node[font=\scriptsize,anchor=east] at (6.18,0.85) {\(u_1\)};
\node[font=\scriptsize,anchor=east] at (6.18,1) {\(1\)};
\end{tikzpicture}
\caption{The two partitions and their common refinement for the case \(L_n=Q=3\) in Figure~\ref{fig:boundary-records}. Each refinement cell \(A_p\cap U_q\) corresponds to the right-crossing edge \(\{\ell_p,r_q\}\).}
\label{fig:boundary-common-refinement}
\end{figure}
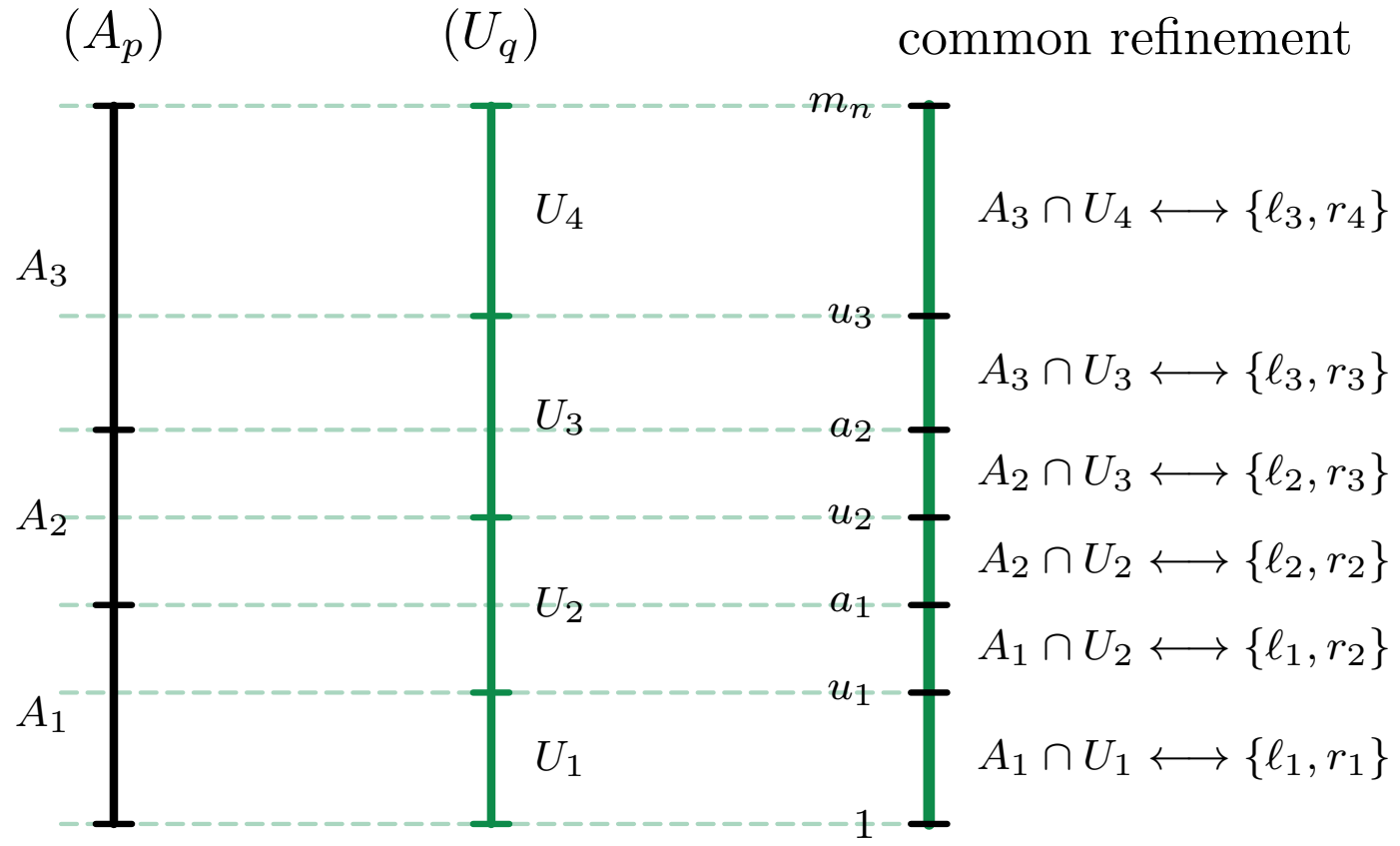

\begin{proof}
Every right-crossing edge has endpoints \(\ell_p,r_q\), since any other endpoint is blocked by an intermediate vertex with a smaller insertion time.
Moreover, only \(1\le q\le Q+1\) can occur: if \(q\ge Q+2\), then \(r_{Q+1}\) lies between \(\ell_p\) and \(r_q\), and \(u_{Q+1}<m_n\le a_p,\) so \(r_{Q+1}\) blocks the edge.
For \(1\le p\le L_n\) and \(1\le q\le Q+1\), the definitions of \(A_p,U_q\), the record construction, and~\eqref{inf-rule} give
\[
A_p\cap U_q\text{ is nondegenerate}
\quad\Longleftrightarrow\quad
\max\{a_p,u_q\}<\min\{a_{p-1},u_{q-1}\},
\]
that is, \(\{\ell_p,r_q\}\in E_\infty\).
The endpoints \(a_0 = u_0 = 1\) handle an empty intervening interval because all insertion times are strictly smaller than \(1\).
Thus the right-crossing edges correspond bijectively to the cells of the common refinement.
\end{proof}

\begin{proof}[Proof of Theorem~\ref{thm:boundary-fluctuations}]
The first partition has \(L_n-1\) interior cut points and the second has \(Q\).
These cut points are distinct almost surely.
By Lemma~\ref{lem:boundary-common-refinement}, every cell corresponds to one right-crossing edge, so
\begin{equation}
\label{eq:right-boundary-decomp}
B_n^+=1+(L_n-1)+Q=L_n+Q=L_n+R_n^+(m_n).
\end{equation}

The relative order of the insertion times in \([n]\) is independent of \(m_n\) and of the insertion times outside \([n]\).
Since \(L_n\) depends only on this relative order, \(L_n\) and \(R_n^+(m_n)\) are independent.

We now prove~\eqref{eq:boundary-as-order}.
Let \(0<s<2\).
By~\eqref{eq:right-boundary-decomp} and the independence described above,
\[
\E\left[s^{B_n^+}\right]
=\E\left[s^{L_n}\right]\E\left[s^{R_n^+(m_n)}\right].
\]
Lemma~\ref{lem:record-generating-function} gives \(\E\left[s^{L_n}\right]=\Gamma(n+s)/(\Gamma(s)\Gamma(n+1))\).
Conditional on \(m_n\), the random variable \(R_n^+(m_n)\) is Poisson with mean \(\log(1/m_n)\), while \(m_n\) has density \(n(1-x)^{n-1}\) on \((0,1)\).
Therefore,
\[
\begin{aligned}
\E\left[s^{R_n^+(m_n)}\right]
&=\E\left[\exp\left\{(s-1)\log(1/m_n)\right\}\right]
=\E\left[m_n^{1-s}\right] \\
&=n\int_0^1 x^{1-s}(1-x)^{n-1}dx
=\frac{\Gamma(n+1)\Gamma(2-s)}{\Gamma(n+2-s)}.
\end{aligned}
\]
Combining these formulas, we obtain
\[
\E\left[s^{B_n^+}\right]
=\frac{\Gamma(2-s)}{\Gamma(s)}
  \frac{\Gamma(n+s)}{\Gamma(n+2-s)}.
\]
For each fixed \(s\in(0,2)\), Fact~\ref{fact:gamma-ratio} gives \(\E\left[s^{B_n^+}\right]\leq C_sn^{2s-2}\) for all sufficiently large \(n\), and increasing \(C_s\) makes the bound valid for all \(n\).
For the upper bound, choose \(s\in(1,2)\).
Since \(s^x\) is increasing, Markov's inequality gives
\[
\P(B_n^+\geq C\log n)
\leq s^{-C\log n}\E\left[s^{B_n^+}\right]
\leq C_sn^{2s-2-C\log s}.
\]
Choosing \(C\) sufficiently large makes the exponent smaller than \(-1\), so these probabilities are summable.
Borel--Cantelli then gives \(B_n^+=O(\log n)\) almost surely.
For the lower bound, choose \(s\in(0,1/2)\).
Now \(s^x\) is decreasing, and Markov's inequality gives
\[
\P(B_n^+\leq c\log n)
\leq s^{-c\log n}\E\left[s^{B_n^+}\right]
\leq C_sn^{2s-2-c\log s}.
\]
Because \(2s-2<-1\) and \(\log s<0\), choosing \(c>0\) sufficiently small again makes the exponent smaller than \(-1\).
Thus the lower tail probabilities are summable, and Borel--Cantelli gives \(B_n^+=\Omega(\log n)\) almost surely.

Consequently, \(B_n^+=\Theta(\log n)\) almost surely.
By symmetry, the same bounds hold for \(B_n^-\).
Hence
\[
X_n=\Theta(\log n)
\qquad\text{a.s}.
\]

We next establish~\eqref{eq:boundary-total-pgf} and prove~\eqref{eq:boundary-gaussian-limit}.
The identity~\eqref{eq:right-boundary-decomp} for right-crossing edges and its left analogue give
\[
X_n=\underbrace{L_n+R_n}_{S_n}
+\underbrace{L_n^-(m_n)+R_n^+(m_n)}_{P_n}.
\]
The minimum location \(\kappa_n\) is uniform on \([n]\), and the two sides have independent relative orders conditional on \(\kappa_n=k\). Hence
\[
S_n\mid(\kappa_n=k)\stackrel d=2+Z_{k-1}^{(1)}+Z_{n-k}^{(2)},
\]
where \(Z_m\) is the number of records in \(m\) i.i.d.\ continuous variables.
Lemma~\ref{lem:record-generating-function} therefore gives
\[
\E\left[s^{S_n}\right]
=\frac{s^2}{n}\sum_{k=1}^nC_{k-1}(s)C_{n-k}(s)
=\frac{s^2}{n}[z^{n-1}](1-z)^{-2s}
=s^2\frac{\Gamma(n-1+2s)}{\Gamma(2s)\Gamma(n+1)}.
\]
The insertion times on the two exterior half lines are independent of each other and of \(m_n\).
Thus the one-sided conditional Poisson law above and its left analogue give, for \(0<x<1\),
\[
L_n^-(m_n),R_n^+(m_n)\mid(m_n=x)
\stackrel{\mathrm{ind}}{\sim}\operatorname{Poisson}\!\bigl(\log(1/x)\bigr).
\]
Consequently,
\[
P_n\mid m_n\sim\operatorname{Poisson}\bigl(2\log(1/m_n)\bigr).
\]
Using the density of \(m_n\) given above, for \(s<3/2\),
\[
\E\left[s^{P_n}\right]
=n\int_0^1x^{2-2s}(1-x)^{n-1}\,dx
=\frac{\Gamma(n+1)\Gamma(3-2s)}{\Gamma(n+3-2s)}.
\]
The independence noted above also shows that \(S_n\) and \(P_n\) are independent.
Multiplying these two generating functions proves~\eqref{eq:boundary-total-pgf}.
Taking the first two derivatives of \(\log\E\left[e^{tX_n}\right]\) at \(t=0\) yields
\[
\E\left[X_n\right]=4H_n,\qquad \operatorname{Var}(X_n)=4H_n+2.
\]

For the central limit theorem, fix \(t\in\R\) and set
\[
s_n=\exp\left\{\frac{t}{2\sqrt{\log n}}\right\}.
\]
For large \(n\), \eqref{eq:boundary-total-pgf} applies.
Since \(s_n\to1\), the prefactor \(s_n^2\Gamma(3-2s_n)/\Gamma(2s_n)\) tends to \(1\).
Moreover, the parameters \(-1+2s_n\) and \(3-2s_n\) remain in a fixed compact interval, so Fact~\ref{fact:gamma-ratio} gives
\[
\frac{\Gamma(n-1+2s_n)}{\Gamma(n+3-2s_n)}
=n^{4s_n-4}(1+o(1)).
\]
Hence
\[
\log\E\left[s_n^{X_n}\right]=(4s_n-4)\log n+o(1),
\]
and therefore
\begin{equation}
\label{eq:boundary-fluctuation-log-mgf}
\log\E\left[\exp\left\{
t\frac{X_n-4\log n}{2\sqrt{\log n}}
\right\}\right]
=-2t\sqrt{\log n}+(4s_n-4)\log n+o(1).
\end{equation}
The Taylor expansion \(e^x=1+x+x^2/2+o(x^2)\) at \(x=0\), with \(x=t/(2\sqrt{\log n})\), gives
\[
(4s_n-4)\log n
=2t\sqrt{\log n}+\frac{t^2}{2}+o(1).
\]
Thus the terms of order \(\sqrt{\log n}\) cancel in~\eqref{eq:boundary-fluctuation-log-mgf}, leaving \(t^2/2+o(1)\).
Consequently, the moment generating functions converge to \(e^{t^2/2}\), and
\[
\frac{|\partial_E[n]|-4\log n}{2\sqrt{\log n}}
\Longrightarrow\mathcal N(0,1).
\]
\end{proof}

\begin{remark}
Since \(|[n]|=n\) and Theorem~\ref{thm:boundary-fluctuations} gives \(|\partial_E[n]|=\Theta(\log n)\) almost surely,
\[
\frac{|\partial_E[n]|}{|[n]|}
=O\left(\frac{\log n}{n}\right)
\longrightarrow0
\qquad\text{almost surely}.
\]
Thus the intervals \([n]\) form a \emph{Følner sequence}\footnote{A sequence of nonempty finite vertex sets \((A_j)\) is a Følner sequence in the edge-boundary sense if \(|\partial_E A_j|/|A_j|\to0\), where \(\partial_E A_j\) consists of edges with exactly one endpoint in \(A_j\).} for \(G_\infty\) almost surely.
\end{remark}

\section{The \texorpdfstring{\(K\)}{K}-Nearest-Predecessor Graph}

The two-sided nearest-predecessor graph connects each new vertex to its nearest predecessor on each side.
A natural extension is to connect each new vertex to its \(K\) nearest predecessors, without prescribing how many lie on either side.
This \(K\)-nearest-predecessor model allows the local density to vary with \(K\) and tests which limiting properties persist without the left--right constraint.

For distinct vertices \(u,v\) in the graph under consideration, we say that \(u\) is \emph{available} to \(v\) if it was inserted before \(v\), that is, \(T_u<T_v\) (equivalently, \(\tau_n(u)<\tau_n(v)\) in the finite model).
We say that \(v\) \emph{selects} \(u\) if \(u\in P_{n,\tau}^{(K)}(v)\) in the finite model or \(u\in P_\infty^{(K)}(v)\) in the infinite model (Definitions~\ref{def:finite-K-nearest-predecessor} and~\ref{def:infinite-K-nearest-predecessor}).
Thus, a later vertex is not available for selection, and an edge joins two vertices if and only if one selects the other.
An available vertex need not be selected, since \(K\) other available vertices may precede it in the order \(\prec_v\).

Because predecessors outside a finite interval may displace vertices selected within it, projective consistency generally fails.
Instead, we show that all edges incident to a fixed vertex, and consequently every fixed-radius rooted ball, are determined by the insertion times in an almost surely finite random interval (Lemmas~\ref{lem:directionless-finite-interval-incident-edges} and~\ref{lem:directionless-finite-interval-rooted-balls}).

\subsection{Nonplanarity}
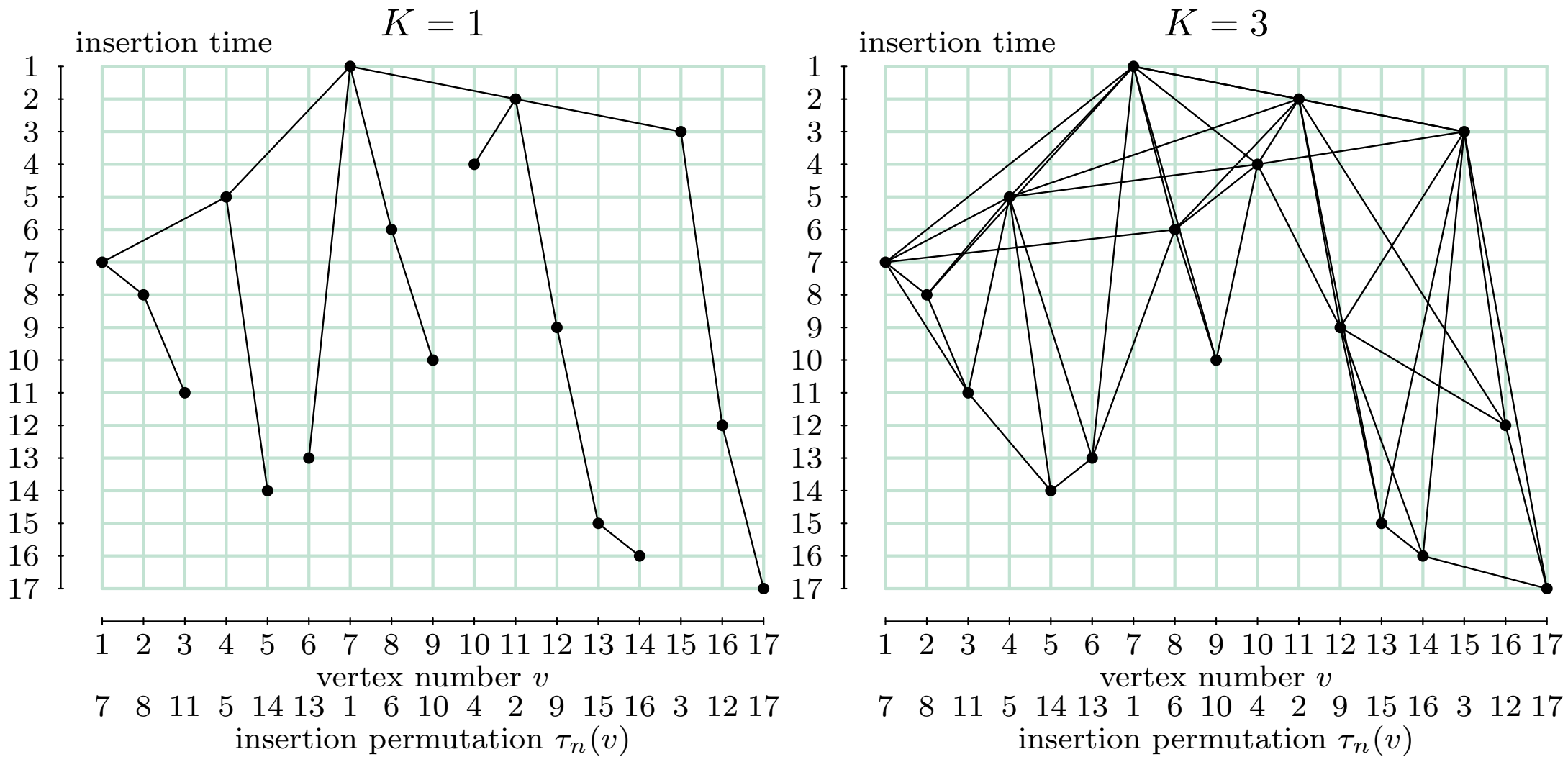
\begin{figure}[ht]
\centering
\begin{tikzpicture}[x=0.38cm,y=-0.30cm,line cap=round,line join=round]

\begin{scope}
\node[anchor=south west,font=\scriptsize] at (0,1) {insertion time};
\node[anchor=south,font=\small] at (9,0.5) {\(K=1\)};

\draw (0,1) -- (0,17);
\foreach \t in {1,...,17} {
    \draw (-0.06,\t) -- (0.06,\t);
    \node[anchor=east,font=\scriptsize] at (-0.16,\t) {\t};
}

\draw (1,18) -- (17,18);
\foreach \v in {1,...,17} {
    \draw (\v,17.9) -- (\v,18.1);
    \node[anchor=base,font=\scriptsize] at (\v,19) {\v};
}
\node[anchor=base,font=\scriptsize] at (9,20) {vertex number \(v\)};

\foreach \v/\t in {1/7,2/8,3/11,4/5,5/14,6/13,7/1,8/6,9/10,10/4,11/2,12/9,13/15,14/16,15/3,16/12,17/17} {
    \node[anchor=base,font=\scriptsize] at (\v,20.9) {\t};
}
\node[anchor=base,font=\scriptsize] at (9,21.9) {insertion permutation \(\tau_n(v)\)};

\foreach \v in {1,...,17} {
    \draw[babyblue!25,line width=0.8pt] (\v,1) -- (\v,17);
}
\foreach \t in {1,...,17} {
    \draw[babyblue!25,line width=0.8pt] (1,\t) -- (17,\t);
}

\foreach \xa/\ya/\xb/\yb in {
1/7/4/5,
1/7/2/8,
2/8/3/11,
4/5/7/1,
4/5/5/14,
6/13/7/1,
7/1/8/6,
8/6/9/10,
10/4/11/2,
7/1/11/2,
11/2/12/9,
12/9/13/15,
13/15/14/16,
11/2/15/3,
15/3/16/12,
16/12/17/17
} {
    \draw[line width=0.45pt] (\xa,\ya) -- (\xb,\yb);
}

\foreach \v/\t in {1/7,2/8,3/11,4/5,5/14,6/13,7/1,8/6,9/10,10/4,11/2,12/9,13/15,14/16,15/3,16/12,17/17} {
    \fill (\v,\t) circle (1.5pt);
}
\end{scope}

\begin{scope}[xshift=7.2cm]
\node[anchor=south west,font=\scriptsize] at (0,1) {insertion time};
\node[anchor=south,font=\small] at (9,0.5) {\(K=3\)};

\draw (0,1) -- (0,17);
\foreach \t in {1,...,17} {
    \draw (-0.06,\t) -- (0.06,\t);
    \node[anchor=east,font=\scriptsize] at (-0.16,\t) {\t};
}

\draw (1,18) -- (17,18);
\foreach \v in {1,...,17} {
    \draw (\v,17.9) -- (\v,18.1);
    \node[anchor=base,font=\scriptsize] at (\v,19) {\v};
}
\node[anchor=base,font=\scriptsize] at (9,20) {vertex number \(v\)};

\foreach \v/\t in {1/7,2/8,3/11,4/5,5/14,6/13,7/1,8/6,9/10,10/4,11/2,12/9,13/15,14/16,15/3,16/12,17/17} {
    \node[anchor=base,font=\scriptsize] at (\v,20.9) {\t};
}
\node[anchor=base,font=\scriptsize] at (9,21.9) {insertion permutation \(\tau_n(v)\)};

\foreach \v in {1,...,17} {
    \draw[babyblue!25,line width=0.8pt] (\v,1) -- (\v,17);
}
\foreach \t in {1,...,17} {
    \draw[babyblue!25,line width=0.8pt] (1,\t) -- (17,\t);
}

\foreach \xa/\ya/\xb/\yb in {
1/7/4/5,
1/7/7/1,
1/7/8/6,
1/7/2/8,
2/8/4/5,
2/8/7/1,
2/8/3/11,
3/11/4/5,
1/7/3/11,
4/5/7/1,
4/5/10/4,
4/5/11/2,
4/5/5/14,
5/14/6/13,
3/11/5/14,
6/13/7/1,
4/5/6/13,
6/13/8/6,
7/1/8/6,
8/6/10/4,
8/6/11/2,
8/6/9/10,
9/10/10/4,
7/1/9/10,
10/4/11/2,
7/1/10/4,
10/4/15/3,
7/1/11/2,
11/2/12/9,
10/4/12/9,
12/9/15/3,
12/9/13/15,
11/2/13/15,
13/15/15/3,
13/15/14/16,
14/16/15/3,
12/9/14/16,
11/2/15/3,
7/1/15/3,
15/3/16/12,
12/9/16/12,
11/2/16/12,
16/12/17/17,
15/3/17/17,
14/16/17/17
} {
    \draw[line width=0.45pt] (\xa,\ya) -- (\xb,\yb);
}

\foreach \v/\t in {1/7,2/8,3/11,4/5,5/14,6/13,7/1,8/6,9/10,10/4,11/2,12/9,13/15,14/16,15/3,16/12,17/17} {
    \fill (\v,\t) circle (1.5pt);
}
\end{scope}
\end{tikzpicture}
\caption{Realizations of finite \(K\)-nearest-predecessor graphs with \(K=1\) (left) and \(K=3\) (right) on \(n=17\) vertices, using the same insertion permutation as Figure~\ref{fig:finite-graph-example}.}
\label{fig:finite-K3-example}
\end{figure}

The finite model and its tie-breaking order were defined in Definition~\ref{def:finite-K-nearest-predecessor}; see Figure~\ref{fig:finite-K3-example} for examples with \(K=1\) and \(K=3\).

Recall the infinite model from Definition~\ref{def:infinite-K-nearest-predecessor}.

\begin{proposition}
\label{prop:directionless-nonplanarity}
For every \(K\geq2\), the graph \(G_\infty^{(K)}\) is nonplanar almost surely.
\end{proposition}

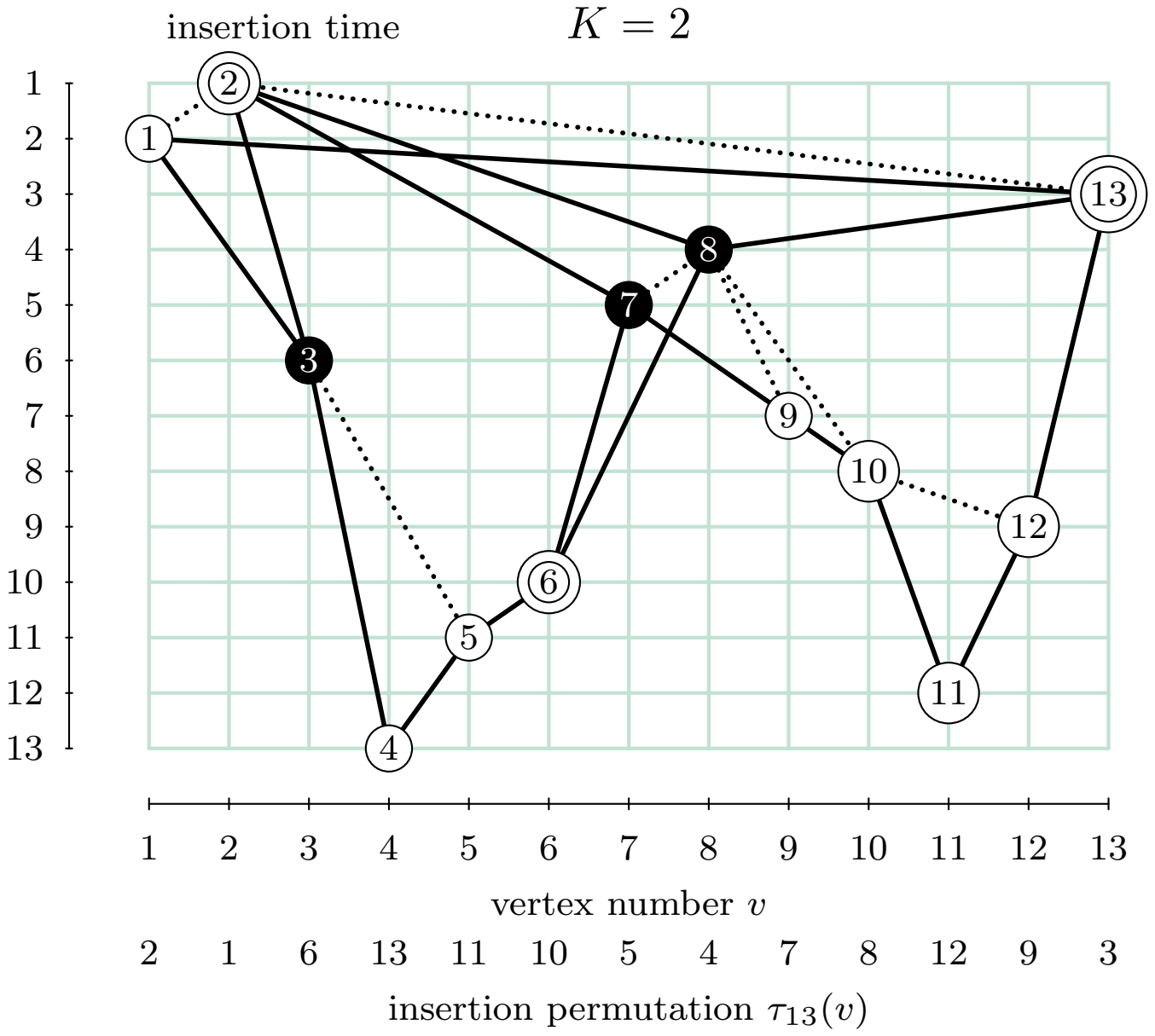
\begin{figure}[ht]
\centering
\begin{tikzpicture}[x=0.62cm,y=-0.43cm,line cap=round,line join=round]
\node[anchor=south west,font=\scriptsize] at (1,0.5) {insertion time};
\node[anchor=south,font=\small] at (7,0.5) {\(K=2\)};

\draw (0,1) -- (0,13);
\foreach \t in {1,...,13} {
    \draw (-0.04,\t) -- (0.04,\t);
    \node[anchor=east,font=\scriptsize] at (-0.10,\t) {\t};
}

\draw (1,14) -- (13,14);
\foreach \v in {1,...,13} {
    \draw (\v,13.9) -- (\v,14.1);
    \node[anchor=base,font=\scriptsize] at (\v,15) {\v};
}
\node[anchor=base,font=\scriptsize] at (7,16) {vertex number \(v\)};

\foreach \v/\t in {1/2,2/1,3/6,4/13,5/11,6/10,7/5,8/4,9/7,10/8,11/12,12/9,13/3} {
    \node[anchor=base,font=\scriptsize] at (\v,16.9) {\t};
}
\node[anchor=base,font=\scriptsize] at (7,17.9)
{insertion permutation \(\tau_{13}(v)\)};

\foreach \v in {1,...,13} {
    \draw[babyblue!25,line width=0.8pt] (\v,1) -- (\v,13);
}
\foreach \t in {1,...,13} {
    \draw[babyblue!25,line width=0.8pt] (1,\t) -- (13,\t);
}

\foreach \xa/\ya/\xb/\yb in {
1/2/2/1,
2/1/13/3,
7/5/8/4,
8/4/9/7,
8/4/10/8,
10/8/12/9,
3/6/5/11
} {
    \draw[black,line width=1.05pt,dash pattern=on 0pt off 2.4pt,line cap=round] (\xa,\ya) -- (\xb,\yb);
}

\foreach \xa/\ya/\xb/\yb in {
3/6/2/1,
3/6/4/13,
4/13/5/11,
5/11/6/10,
3/6/1/2,
1/2/13/3,
7/5/2/1,
7/5/6/10,
7/5/9/7,
9/7/10/8,
10/8/11/12,
11/12/12/9,
12/9/13/3,
8/4/2/1,
8/4/6/10,
8/4/13/3
} {
    \draw[line width=1.05pt] (\xa,\ya) -- (\xb,\yb);
}

\foreach \v/\t in {1/2,4/13,5/11,9/7,10/8,11/12,12/9} {
    \node[draw,circle,fill=white,inner sep=1.25pt,font=\scriptsize]
    at (\v,\t) {\v};
}

\foreach \v/\t in {3/6,7/5,8/4} {
    \node[draw,circle,fill=black,text=white,inner sep=1.25pt,font=\scriptsize]
    at (\v,\t) {\v};
}
\foreach \v/\t in {2/1,6/10,13/3} {
    \node[draw,circle,fill=white,inner sep=2.5pt,font=\scriptsize]
    at (\v,\t) {\phantom{\v}};
    \node[draw,circle,fill=white,inner sep=0.8pt,font=\scriptsize]
    at (\v,\t) {\v};
}
\end{tikzpicture}
\caption{A nonplanar realization for \(K=2\).
The solid edges form a subdivision of \(K_{3,3}\) with parts \(\{3,7,8\}\) (black) and \(\{2,6,13\}\) (double-circled); the dotted edges are the remaining edges.}
\label{fig:directionless-nonplanarity}
\end{figure}

\begin{proof}
We first consider \(K=2\).
In Figure~\ref{fig:directionless-nonplanarity}, the solid edges give nine connections, one between each vertex in \(\{3,7,8\}\) and each vertex in \(\{2,6,13\}\).
Each is a single edge or a path through white vertices, and no white vertex lies on two paths.
Thus the solid subgraph is a subdivision of \(K_{3,3}\) and is nonplanar by Kuratowski's theorem.

To force this configuration in the infinite graph, consider \(24\) consecutive vertices, numbered \(1,\ldots,24\) from the left.
Require vertices \(1,\ldots,13\) to have the relative insertion order shown in the figure and to precede vertices \(14,\ldots,24\).
This event has positive probability because it imposes only finitely many strict inequalities among continuous i.i.d.\ insertion times.

The relative order of vertices \(1,\ldots,13\) determines all required solid edges except for the choices made by vertex \(13\).
The possible competitors for vertex \(13\) are vertices \(14,\ldots,24\), all of which are inserted later, while vertices outside \(\{1,\ldots,24\}\) come after its chosen predecessors in the order \(\prec_{13}\).
Hence, the solid subdivision is present regardless of the insertion times outside the block.

Apply the construction to disjoint intervals of length \(24\).
The corresponding events are independent and have the same positive probability, so the second Borel--Cantelli lemma yields such a subdivision almost surely.
Therefore, \(G_\infty^{(2)}\) is nonplanar; since \(E_\infty^{(2)}\subseteq E_\infty^{(K)}\), the same holds for every \(K\ge2\).
\end{proof}

\begin{remark}
Finding the 13-vertex configuration in Figure~\ref{fig:directionless-nonplanarity} is not by itself enough in the infinite model: predecessors outside the interval may displace vertices selected within it.
For the same reason, unlike the two-sided nearest-predecessor model, we generally do not have
\[
G_n^{(K)} = G_\infty^{(K)}[[n]].
\]
We therefore use the \emph{blocker intervals} \(I_v\) introduced in Lemma~\ref{lem:directionless-finite-interval-incident-edges} to determine all edges incident to a fixed vertex; Lemma~\ref{lem:directionless-finite-interval-rooted-balls} extends this determination to rooted balls by induction on graph distance.
\end{remark}

\subsection{The Quenched BS Limit}
For \(v\in\Z\), define the \emph{one-sided blocker distances}
\[
\begin{aligned}
\beta_K^-(v)
&=\min\left\{r\geq1:\sum_{j=1}^r \mathbf 1\{T_{v-j}<T_v\}\geq K\right\},\\
\beta_K^+(v)
&=\min\left\{r\geq1:\sum_{j=1}^r \mathbf 1\{T_{v+j}<T_v\}\geq K\right\}.
\end{aligned}
\]
Conditional on \(T_v=t>0\), each half line contains infinitely many vertices with insertion times smaller than \(t\).
Hence, almost surely, for every \(v\in\Z\),
\[
\beta_K^-(v)<\infty, \qquad \beta_K^+(v)<\infty.
\]
Set
\[
Q_K(v) = \max\{\beta_K^-(v),\beta_K^+(v)\}, \qquad N_G(v) = \{u\in V(G):\{u,v\}\in E(G)\}.
\]

\begin{lemma}
\label{lem:directionless-finite-interval-incident-edges}
For \(v\in\Z\), set
\[
I_v=[v-Q_K(v),v+Q_K(v)].
\]
Almost surely,
\[
N_{G_\infty^{(K)}}(v)\subseteq I_v.
\]
Moreover, for every \(n\geq1\) and \(v\in[n]\), if \(I_v\subseteq[n]\), then
\begin{equation}
\label{eq:directionless-blocker-interval-neighborhood-consistency}
N_{G_n^{(K)}}(v)=N_{G_\infty^{(K)}}(v).
\end{equation}
Thus all edges incident to \(v\) are determined by the insertion times in \(I_v\).
\end{lemma}

\begin{proof}
If \(u\notin I_v\), then, by the definition of the blocker distances, the interval \(I_v\) contains at least \(K\) \emph{blockers} lying between \(v\) and \(u\), all with insertion times smaller than \(T_v\).
They precede \(u\) in \(\prec_v\), so \(v\) cannot select \(u\) when \(T_u<T_v\); when \(T_v<T_u\), \(u\) is not available to \(v\).
For the reverse choice, if \(T_v>T_u\), then \(v\) is not available to \(u\); if \(T_v<T_u\), the \(K\) blockers have insertion times smaller than \(T_v<T_u\) and precede \(v\) in \(\prec_u\), so \(u\) cannot select \(v\).
Hence, \(u\) and \(v\) are not adjacent in \(G_\infty^{(K)}\).
If \(u,v\in[n]\), all the blockers between them also belong to \([n]\), so \(u\) and \(v\) are not adjacent in \(G_n^{(K)}\) either.
Consequently,
\begin{equation}
\label{eq:directionless-blocker-interval-exclusion}
|u-v|>Q_K(v)
\quad\Longrightarrow\quad
u\notin N_{G_\infty^{(K)}}(v),
\quad
u\notin N_{G_n^{(K)}}(v)\ \text{if }u,v\in[n].
\end{equation}

Now suppose that \(I_v\subseteq[n]\).
By the definition of \(Q_K(v)\), the first \(K\) vertices with smaller insertion times in the order \(\prec_v\) all lie in \(I_v\).
Therefore,
\begin{equation}
\label{eq:directionless-blocker-interval-outgoing}
P_{n,\tau}^{(K)}(v)=P_\infty^{(K)}(v).
\end{equation}

Fix \(u\in I_v\setminus\{v\}\).
If \(T_u<T_v\), then \(v\) is unavailable to \(u\) in both models.
Suppose that \(T_v<T_u\), and consider the \(K\) blockers on the same side of \(v\) as \(u\).
If the farthest blocker is at coordinate distance less than \(2|u-v|\) from \(v\), then all \(K\) blockers precede \(v\) in \(\prec_u\).
Their insertion times are smaller than \(T_v<T_u\), so \(u\) cannot select \(v\) in either model.
Otherwise, the farthest blocker is at distance at least \(2|u-v|\) from \(v\).
Since this blocker lies in \(I_v\), we have \(2|u-v|\leq Q_K(v)\).
Moreover, \(w\prec_u v\) implies \(|w-u|\leq|v-u|\), and hence
\[
|w-v|\leq |w-u|+|u-v|\leq2|u-v|\leq Q_K(v).
\]
Thus every candidate that can affect whether \(u\) selects \(v\) belongs to \(I_v\), so the finite and infinite decisions agree.
We have therefore shown that
\begin{equation}
\label{eq:directionless-blocker-interval-incoming}
v\in P_{n,\tau}^{(K)}(u)
\quad\Longleftrightarrow\quad
v\in P_\infty^{(K)}(u)
\qquad (u\in I_v\setminus\{v\}).
\end{equation}
Combining~\eqref{eq:directionless-blocker-interval-exclusion}, \eqref{eq:directionless-blocker-interval-outgoing}, and~\eqref{eq:directionless-blocker-interval-incoming} proves~\eqref{eq:directionless-blocker-interval-neighborhood-consistency}.
\end{proof}

\begin{corollary}
\label{cor:directionless-local-finiteness}
The graph \(G_\infty^{(K)}\) is locally finite almost surely.
\end{corollary}

\noindent The proof is given in Appendix~\ref{app:proof-directionless-local-finiteness}.

For \(r\in\N_0\) and \(v\in\Z\), define
\[
\mathcal R_r(v) = \max_{x\in B_r(G_\infty^{(K)},v)} \left( |x-v|+Q_K(x) \right).
\]
By Corollary~\ref{cor:directionless-local-finiteness} and the finiteness of the \(Q_K(x)\)'s, \(\mathcal R_r(v)<\infty\) almost surely; its coordinate interval contains the determining intervals of every vertex in the rooted ball.

\begin{lemma}
\label{lem:directionless-finite-interval-rooted-balls}
Almost surely, for every \(n\geq1\), \(v\in[n]\), and \(r\in\N_0\), if
\[
[v-\mathcal R_r(v),v+\mathcal R_r(v)] \subseteq[n],
\]
then
\[
B_r(G_n^{(K)},v) = B_r(G_\infty^{(K)},v).
\]
\end{lemma}

\noindent The proof is given in Appendix~\ref{app:proof-directionless-finite-interval-rooted-balls}.

Let \(\mathcal C_\infty^{(K)}(0)\) be the component of \(0\) in \(G_\infty^{(K)}\), and define
\[
\nu_n^{(K)}=\frac1n\sum_{v=1}^n\delta_{[G_n^{(K)},v]},
\qquad
\mu_K=\mathcal L\left([\mathcal C_\infty^{(K)}(0),0]\right).
\]
Theorem~\ref{thm:directionless-one-endedness} shows that \(\mathcal C_\infty^{(K)}(0)=G_\infty^{(K)}\) almost surely for \(K \ge 1\).

\begin{proof}[Proof of Theorem~\ref{thm:directionless-quenched-BS}]
The proof is similar to the proof of Theorem~\ref{thm:quenched-BS}.

Fix \(r\in\N_0\) and a finite rooted graph \(H\).
For the insertion time shift \(\theta\), set
\[
X_v= \mathbf 1\{B_r(G_\infty^{(K)},v)\cong H\}=X_0\circ\theta^v,
\qquad
Y_{n,v}=\mathbf 1\{B_r(G_n^{(K)},v)\cong H\}.
\]
Birkhoff's theorem gives \(n^{-1}\sum_{v=1}^nX_v\to\E\left[X_0\right]\) almost surely.
For \(h\in\N\), Lemma~\ref{lem:directionless-finite-interval-rooted-balls} gives
\[
\left| \frac1n\sum_{v=1}^nY_{n,v} - \frac1n\sum_{v=1}^nX_v \right| \leq \frac{2h}{n} + \frac1n \sum_{v=1}^n \mathbf 1\{\mathcal R_r(v)>h\}.
\]
Since \(\mathcal R_r(v)=\mathcal R_r(0)\circ\theta^v\), applying Birkhoff to the last term and then letting \(h\to\infty\) shows that the two averages differ by \(o(1)\) almost surely.
Thus, the desired rooted ball probabilities converge; a countable intersection over \((r,H)\) proves the claim by Lemma~\ref{lem:cylinder-events-determine-BS-convergence}.
\end{proof}

\subsection{One-Endedness and Recurrence Without Planarity}
\label{sec:directionless-recurrence}

For \(d\geq1\), write
\[
p_K(d) = \P\bigl(\{0,d\}\in E_\infty^{(K)}\bigr).
\]

\begin{lemma}
\label{lem:directionless-edge-probability}
For every \(d\ge1\),
\[
p_K(d)\leq \frac{K^2}{d^2}.
\]
If, in addition, \(2d-1\geq K\), then
\[
p_K(d)=\frac{K(K+1)}{4d^2-1}.
\]
\end{lemma}

\begin{proof}
Suppose that a candidate \(y\) is preceded by \(m\) vertices in the distance order of a source \(x\).
Conditional on \(T_x=t\), the source selects \(y\) when \(T_y<t\) and at most \(K-1\) of those \(m\) preceding vertices have insertion times below \(t\).
Setting \(q=\min\{K,m+1\}\), we obtain
\[
\P\bigl(y\in P_\infty^{(K)}(x)\bigr) = \sum_{r=0}^{q-1} \binom mr \int_0^1 t^{r+1}(1-t)^{m-r}\,dt
=\frac{q(q+1)}{2(m+1)(m+2)}.
\]
In the order \(\prec_d\), the vertex \(0\) is preceded by \(2d-2\) candidates; in the order \(\prec_0\), the vertex \(d\) is preceded by \(2d-1\) candidates.
The two possible orientations of the edge are disjoint, since they require opposite orderings of \(T_0\) and \(T_d\).
When \(2d-1\geq K\), adding their probabilities gives
\[
p_K(d) = \frac{K(K+1)}{2(2d-1)(2d)} + \frac{K(K+1)}{2(2d)(2d+1)} = \frac{K(K+1)}{4d^2-1}.
\]
This is at most \(K^2/d^2\).
If \(2d-1<K\), then \(p_K(d)\leq1\leq K^2/d^2\), proving the uniform bound.
\end{proof}

\begin{remark}
The left--right bypass argument does not transfer directly because the directionless rule imposes no left--right allocation: a vertex may select all \(K\) predecessors on the same side.
\end{remark}

\begin{theorem}
\label{thm:directionless-one-endedness}
For every fixed integer \(K\geq1\), almost surely \(G_\infty^{(1)}\) is a one-ended spanning tree of \(G_\infty^{(K)}\).
Consequently, \(G_\infty^{(K)}\) is almost surely connected and one-ended.
\end{theorem}

\begin{proof}
We first analyze the case \(K=1\).
For \(v\in\Z\), write \(p(v)\) for the unique element of \(P_\infty^{(1)}(v)\), and orient the edge \(\{v,p(v)\}\) from \(v\) to \(p(v)\).
Thus every vertex has outdegree one and \(T_{p(v)}<T_v\).

If the underlying undirected graph contained a cycle, the vertex on the cycle with the largest insertion time would have both incident cycle edges oriented away from it, contradicting its outdegree one.
Hence \(G_\infty^{(1)}\) is a forest.
Moreover, for every \(v\), we obtain an infinite ancestral ray:
\[
v, \ p(v), \ p(p(v)), \ \ldots.
\]

We next show that all descendant subtrees are finite.
For \(m\geq0\), set
\[
Z_m(v) = \#\{u\in\Z:p^m(u)=v\}.
\]
Conditional on a vertex having insertion time \(s\), its parent's insertion time is uniform on \((0,s)\), independently of the exploration history.
Thus, if \(U_0,\ldots,U_m \stackrel{\mathrm{ind}}\sim \operatorname{Unif}(0,1)\), then, for \(0<t<1\),
\begin{equation}
\label{eq:directionless-descendant-generations}
\begin{aligned}
\E\left[Z_m(v)\mid T_v=t\right]
&= \E\left[Z_m(0)\mid T_0=t\right]
= \frac{d}{dt}\E\left[Z_m(0)\mathbf 1_{\{T_0\leq t\}}\right] \\
&= \frac{d}{dt}\P(T_{p^m(0)}\leq t)
= \frac{d}{dt}\P\left(\prod_{j=0}^m U_j\leq t\right)\\
&= \frac{(-\log t)^m}{m!}.
\end{aligned}
\end{equation}
The equalities use translation invariance, the uniform density of \(T_0\), a mass transport count by \(m\)-th ancestors, and the fact that the \(m\)-th ancestral insertion time has the distribution of a product of \(m+1\) independent uniform random variables.

The events \(\{Z_m(v)>0\}\) decrease with \(m\). Then, \eqref{eq:directionless-descendant-generations} gives
\[
\P\bigl(Z_m(v)>0\mid T_v=t\bigr)
\leq \E\left[Z_m(v)\mid T_v=t\right]
= \frac{(-\log t)^m}{m!}
\longrightarrow0.
\]
Thus,
\(
\P\left(\bigcap_{m\geq0}\{Z_m(v)>0\}\middle|T_v=t\right)
=\lim_{m\to\infty}\P\bigl(Z_m(v)>0\mid T_v=t\bigr)=0.
\)
For each fixed vertex this holds almost surely; by countability, it holds for every vertex on a single event of probability one.
Since each generation is finite by Corollary~\ref{cor:directionless-local-finiteness}, every descendant subtree is finite.
Ancestral rays starting in the same tree component must coalesce; otherwise the finite path connecting two such rays would contain a vertex with two outgoing edges.
Moreover, because every descendant subtree is finite, every infinite ray eventually follows the ancestral direction.
Thus, all rays in a component represent the same end.

If two edges \(\{a,b\}\) and \(\{c,d\}\) crossed, with \(a<c<b<d\), the first edge would imply \(T_c>T_b\), whereas the second would imply \(T_b>T_c\), a contradiction.
Thus the edges are noncrossing.

Since every component contains an infinite ancestral ray, noncrossing implies that the vertex set of every component is an interval in \(\Z\).
To see this explicitly, suppose \(x<y<z\), with \(x,z\) in one component and \(y\) in another.
A path from \(x\) to \(z\) contains an edge \(\{a,b\}\) with \(a<y<b\).
The component of \(y\) cannot leave the finite interval \((a,b)\): its first edge leaving that interval would either cross \(\{a,b\}\) or meet the component containing \(a\) and \(b\).
This would make the component of \(y\) finite, contradicting its infinite ancestral ray.

There are therefore at most two components.
If there were two, their vertex sets would be
\[
(-\infty,j]\qquad\text{and}\qquad [j+1,\infty)
\]
for a unique \(j\in\Z\).
Let \(A_j\) be the event that these are the two components; the events \((A_j)_{j\in\Z}\) are disjoint and \(\P(A_j)=\P(A_0)\).
Thus, for every \(N\in\N_0\),
\[
1\geq \P\left(\bigcup_{j=-N}^{N}A_j\right)
=\sum_{j=-N}^{N}\P(A_j)
=(2N+1)\P(A_0).
\]
This forces \(\P(A_j)=\P(A_0)=0\) for every \(j\in\Z\).
Countable additivity then gives
\[
\P\left(\bigcup_{j\in\Z}A_j\right)=0.
\]
Then the forest is connected almost surely, hence a tree.

Finally, by definition,
\[
G_\infty^{(1)} \subseteq G_\infty^{(2)} \subseteq \cdots.
\]
Since all these graphs have vertex set \(\Z\), \(G_\infty^{(1)}\) is a one-ended spanning tree of \(G_\infty^{(K)}\) for every \(K\geq1\).

Fix \(K\geq1\) and a finite set \(S\subset\Z\).
Since \(G_\infty^{(1)}\) is locally finite, deleting \(S\) leaves finitely many tree components, exactly one of which is infinite.
The components of \(G_\infty^{(K)}\setminus S\) are obtained by merging these tree components, so exactly one is infinite.
Hence, \(G_\infty^{(K)}\) is connected and one-ended.
\end{proof}

The \(K\)-nearest-predecessor limit is generally nonplanar, so recurrence cannot be deduced from a planar-limit theorem.
\begin{theorem}[{B\"aumler~\cite[Theorem~1.3 and Corollary~3.1]{Baumler2023}}]
\label{thm:baumler-quadratic-kernel}
Let \(\mathcal H\) be a random simple graph on \(\Z\).\footnote{The original corollary is stated for the complete graph on \(\Z\), viewed as an ambient weighted graph in which every pair receives a conductance.
Here we take \(\omega(\{x,y\})=\mathbf 1_{\{\{x,y\}\in E(\mathcal H)\}}\); thus an absent edge has conductance \(0\), and \(\E\left[\omega(\{x,y\})\right]=\P(\{x,y\}\in E(\mathcal H))\leq C/|x-y|^2\).} If there is a constant \(C<\infty\) such that
\[
\P\bigl(\{x,y\}\in E(\mathcal H)\bigr) \leq \frac{C}{|x-y|^2} \qquad (x\neq y),
\]
then every connected component of \(\mathcal H\) is recurrent almost surely.
\end{theorem}

\begin{proposition}
\label{prop:directionless-recurrence}
For every fixed \(K\geq1\), every connected component of \(G_\infty^{(K)}\) is recurrent almost surely.
\end{proposition}

\begin{proof}
Together, translation invariance, Lemma~\ref{lem:directionless-edge-probability}, and Theorem~\ref{thm:baumler-quadratic-kernel} prove the claim.
\end{proof}

\subsection{Degree Tail Conjecture}
The following conjecture predicts geometric decay of the point probabilities of the limiting degree distribution for every fixed \(K\geq1\).
                      
\begin{conjecture}
\label{conj:directionless-K-degree-tail}
For fixed \(K\geq1\), let
\[
D^{(K)}:=\deg_{G_\infty^{(K)}}(0).
\]
There exist \(B_K>0\) and \(\rho_K\in(0,1)\) such that
\[
\P\!\left(D^{(K)}=K+m\right)\sim B_K\rho_K^m
\qquad (m\to\infty).
\]
\end{conjecture}

This conjecture, whose proof would require control of the dependence between later selections, is essentially motivated by Theorem~\ref{thm:degree-distribution-limit-graph}.
Our intuition is as follows.
The root itself selects exactly \(K\) predecessors, so any degree above \(K\) is caused by later vertices that select the root.
The favorable order is that the root arrives before the selecting vertex, while all but at most \(K-1\) of the candidates ahead of the root in that vertex's distance order have yet to arrive when the selecting vertex is inserted.
We therefore condition on a very early root time \(T_0=t\) and fix a potential selecting vertex at a large coordinate distance \(d\) well below \(1/t\).
At this distance, \(2d-2\) candidates lie strictly closer to that vertex than the root.
The potential selecting vertex and its competing candidates retain i.i.d.\ insertion times under this conditioning.
Its chance of arriving before all but at most \(K-1\) competitors is therefore roughly inversely proportional to their number, and hence to distance.
Summing these conditional selection probabilities over the relevant distances suggests that the conditional expected number of later neighbors is of order \(\log(1/t)\).
Since \(-\log T_0\) is exponentially distributed, this suggests exponential degree decay.

\appendix
\addtocontents{toc}{\protect\setcounter{tocdepth}{1}}
\section{Deferred Proofs}
\label{app:deferred-proofs}

\subsection{Proof of Lemma~\ref{lem:cylinder-events-determine-BS-convergence}}
\label{app:proof-cylinder-events}

\begin{proof}
First observe that the cylinder events form a countable basis for the local topology.
Countability follows because there are only countably many finite rooted graphs up to isomorphism.
Every cylinder is both closed and open.
For example, if \((G,o)\in C_{r,H}\), then
\[
d_{\mathrm{loc}}\bigl((G,o),(G',o')\bigr) < \frac1{r+1} \quad\Longrightarrow\quad (G',o')\in C_{r,H}.
\]
The complement of \(C_{r,H}\) is a union of other open radius-\(r\) cylinders, hence is open.
Therefore, \(C_{r,H}\) is also closed.

Let \(O\subseteq\mG_*\) be open and let \((G,o)\in O\).
Since \(O\) is open, there exists \(\varepsilon>0\) such that
\[
\{(G',o'):d_{\mathrm{loc}}((G,o),(G',o'))<\varepsilon\} \subseteq O.
\]
Choose \(r\in\N_0\) such that \(1/(1+r)<\varepsilon\).
If \((G',o')\in C_{r,B_r(G,o)}\), then
\[
d_{\mathrm{loc}}((G,o),(G',o')) \le \frac1{1+r} < \varepsilon.
\]
Therefore,
\[
(G,o)\in C_{r,B_r(G,o)}\subseteq O.
\]
Thus, the cylinders form a basis.

Finite intersections of cylinder events are again either empty or cylinder events.
Indeed, if \(r\le s\),
\[
C_{r,H}\cap C_{s,H'} =
\begin{cases}
C_{s,H'}, & \text{if } B_r(H')\cong H,\\
\emptyset, & \text{otherwise.}
\end{cases}
\]
Now let \(O\subseteq\mG_*\) be open.
Since the cylinder events form a countable basis, there exist cylinder events \(C_1,C_2,\ldots\) such that
\[
O=\bigcup_{j=1}^\infty C_j.
\]
For \(m\ge1\) and each nonempty \(I\subseteq\{1,\ldots,m\}\), define
\[
O_m=\bigcup_{j=1}^m C_j, \qquad C_I=\bigcap_{j\in I}C_j.
\]
By inclusion--exclusion,
\[
\nu_n(O_m) = \sum_{\emptyset\neq I\subseteq\{1,\ldots,m\}} (-1)^{|I|+1}\nu_n(C_I).
\]
Each \(C_I\) is either empty or again a cylinder event.
Hence, by~\eqref{eq:cylinder-convergence}, \(\nu_n(C_I)\longrightarrow \mu(C_I)\).
Since this is a finite sum, we may pass to the limit term by term:
\[
\lim_{n\to\infty}\nu_n(O_m) = \sum_{\emptyset\neq I\subseteq\{1,\ldots,m\}} (-1)^{|I|+1}\mu(C_I) = \mu\left(\bigcup_{j=1}^m C_j\right) = \mu(O_m).
\]
Therefore, for every fixed \(m\),
\[
\nu_n(O_m)\longrightarrow \mu(O_m).
\]
Since \(O_m\subseteq O\), we have
\[
\liminf_{n\to\infty}\nu_n(O) \ge \liminf_{n\to\infty}\nu_n(O_m) = \mu(O_m).
\]
Letting \(m\to\infty\), and using \(O_m\uparrow O\), gives
\[
\liminf_{n\to\infty}\nu_n(O) \ge \mu(O).
\]
By the Portmanteau theorem, \(\nu_n\Longrightarrow \mu\).
\end{proof}

\subsection{Proof of Lemma~\ref{lem:record-generating-function}}
\label{app:proof-record-generating-function}

\begin{proof}
The case \(m=0\) is immediate.
For \(m\ge1\), let \(\tau\) be a uniform permutation on \(\{1,\ldots,m\}\) and set \(I_i=\mathbf 1\{\tau(i)=\min_{j\le i}\tau(j)\}\).
The successive ranks of \(\tau(i)\) among \(\tau(1),\ldots,\tau(i)\) are independent and uniform on \(\{1,\ldots,i\}\); hence the \(I_i\) are independent with \(\P(I_i=1)=1/i\).
Since \(R_m=\sum_{i=1}^m I_i\),
\[
C_m(z)=\prod_{i=1}^m\left(1-\frac1i+\frac zi\right)
=\frac{z(z+1)\cdots(z+m-1)}{m!}.
\]
The gamma function form follows from the rising factorial identity, and reversing \(\tau\) gives the left record formula.
\end{proof}

\subsection{Proof of Lemma~\ref{lem:two-sided-number-records-estimate}}
\label{app:proof-two-sided-number-records}

\begin{proof}
\(L_m\) and \(R_m\) are independent because the two sides use independent relative orders.
By Lemma~\ref{lem:record-generating-function}, for every fixed \(z>1\),
\[
\E\left[z^{S_m}\right] = \E\left[z^{L_m}\right] \E\left[z^{R_m}\right] = C_m(z)^2 = \left( \frac{\Gamma(m+z)} {\Gamma(z)\Gamma(m+1)} \right)^2.
\]
By Fact~\ref{fact:gamma-ratio}, for fixed \(z\),
\begin{equation}
\label{log-E}
\E\left[z^{S_m}\right] = \frac{m^{2(z-1)}}{\Gamma(z)^2}(1+o(1)),\qquad
\log \E\left[z^{S_m}\right] = 2(z-1)\log m+O_z(1).
\end{equation}

We first prove the upper bound.
For any fixed \(z>1\), Markov's inequality gives
\[
\begin{aligned}
\P(S_m\ge u\log m)
&= \P(z^{S_m} \ge z^{u\log m})
\le
z^{-u\log m}\E\left[z^{S_m}\right] \\
&=
\exp\left(
-u(\log z)\log m
+
2(z-1)\log m
+
O_z(1)
\right) \\
&=
m^{-u\log z+2(z-1)+o(1)}.
\end{aligned}
\]
The function \(u\log z-2(z-1)\) is strictly concave and, since \(u>2\), is maximized over \(z>1\) at \(z=u/2\), where its value is \(I(u)\).
Hence
\[
\P(S_m\ge u\log m)
\le
m^{-I(u)+o(1)}.
\]

We now prove the lower bound.
Fix \(u'>u\), and set
\[
z':=\frac{u'}2, \qquad \Lambda_m(z):=\log\E\left[z^{S_m}\right].
\]
Define the \emph{tilted probability measure} \(\Q_m\) by\footnote{Here \(\{S_m\ge u\log m\}\) is the upper tail event, and we are establishing its lower bound in the regime \(u>2\).

Since \(z'=u'/2>1\), this tilt gives more weight to outcomes with larger \(S_m\).
Under \(\P\), the typical size is \(S_m\sim 2\log m\), so \(\{S_m\ge u\log m\}\) is a rare upper tail event.
The tilt is chosen so that under \(\Q_m\), \(S_m\) is centered near \(u'\log m\), making \(S_m\ge u\log m\) typical.
The upper cutoff in \(A_m\) keeps the change of measure factor \(z'^{-S_m}\) from becoming too small when converting back to \(\P\).}
\begin{equation}
\label{measure-Q}
\frac{d\Q_m}{d\P} = \frac{z'^{S_m}}{\E\left[z'^{S_m}\right]}.
\end{equation}
From the exact PGF formula,
\[
\Lambda_m(z) = 2\bigl(\log\Gamma(m+z)-\log\Gamma(m+1)-\log\Gamma(z)\bigr).
\]
Let \(\psi=(\log\Gamma)'=\Gamma'/\Gamma\) be the digamma function and \(\psi_1=\psi'\) the trigamma function.
Then
\[
\Lambda_m'(z)=2\bigl(\psi(m+z)-\psi(z)\bigr),
\qquad
\Lambda_m''(z)=2\bigl(\psi_1(m+z)-\psi_1(z)\bigr).
\]
Hence, uniformly for \(z\) in a compact neighborhood of \(z'\),
\[
\begin{aligned}
\Lambda_m(z) &=2(z-1)\log m-2\log\Gamma(z)+O(m^{-1}),\\
\Lambda_m'(z) &=2\log m-2\psi(z)+O(m^{-1}),\\
\Lambda_m''(z) &=-2\psi_1(z)+O(m^{-1}).
\end{aligned}
\]

Using~\eqref{measure-Q}, for the expectation, we have
\[
\begin{aligned}
\E_{\Q_m}\left[S_m\right] &= \E\left[ S_m\frac{z'^{S_m}}{\E\left[z'^{S_m}\right]} \right] = \frac{ z'\left.\frac{d}{dz}\E\left[z^{S_m}\right]\right|_{z=z'} }{ \E\left[z'^{S_m}\right] } = z'\Lambda_m'(z') = u'\log m+O_{u'}(1).
\end{aligned}
\]
For the variance, first note that
\[
\E_{\Q_m}\left[S_m(S_m-1)\right] = \E\left[ S_m(S_m-1)\frac{z'^{S_m}}{\E\left[z'^{S_m}\right]} \right] = \frac{ z'^2\left.\frac{d^2}{dz^2}\E\left[z^{S_m}\right]\right|_{z=z'} }{ \E\left[z'^{S_m}\right] },
\]
\[
\frac{\left(\E\left[z^{S_m}\right]\right)'}{\E\left[z^{S_m}\right]} = \Lambda_m'(z), \qquad \frac{\left(\E\left[z^{S_m}\right]\right)''}{\E\left[z^{S_m}\right]} = \Lambda_m''(z)+\left(\Lambda_m'(z)\right)^2,
\]
so we have
\[
\E_{\Q_m}\left[S_m^2\right] = \E_{\Q_m}\left[S_m(S_m-1)\right] + \E_{\Q_m}\left[S_m\right] = z'^2\left(\Lambda_m''(z')+\left(\Lambda_m'(z')\right)^2\right) + z'\Lambda_m'(z').
\]
Therefore,
\[
\operatorname{Var}_{\Q_m}(S_m) = \E_{\Q_m}\left[S_m^2\right] - \left(\E_{\Q_m}\left[S_m\right]\right)^2 = z'^2\Lambda_m''(z')+z'\Lambda_m'(z') = u'\log m+O_{u'}(1).
\]

Choose \(A>0\) sufficiently large, and define
\[
A_m:= \left\{ u\log m \le S_m \le u'\log m+A\sqrt{\log m} \right\} \subseteq \{S_m \ge u\log m\}.
\]
Since \(u'>u\), for all sufficiently large \(m\), we have
\[
\E_{\Q_m}\left[S_m\right]-u\log m = (u'-u)\log m+O_{u'}(1) \ge \frac{u'-u}{2}\log m.
\]
Therefore, by Chebyshev's inequality,
\[
\begin{aligned}
\Q_m(S_m<u\log m) &\le \Q_m\left( \left|S_m-\E_{\Q_m}\left[S_m\right]\right| \ge \E_{\Q_m}\left[S_m\right]-u\log m \right)\\
&\le \Q_m\left( \left|S_m-\E_{\Q_m}\left[S_m\right]\right| \ge \frac{u'-u}{2}\log m \right)\\
&\le \frac{\operatorname{Var}_{\Q_m}(S_m)} {\left(\frac{u'-u}{2}\log m\right)^2} = \frac{u'\log m+O_{u'}(1)} {\left(\frac{u'-u}{2}\right)^2(\log m)^2}\\
&\longrightarrow 0.
\end{aligned}
\]
Similarly, for all sufficiently large \(m\),
\[
u'\log m+A\sqrt{\log m}-\E_{\Q_m}\left[S_m\right] = A\sqrt{\log m}+O_{u'}(1) \ge \frac A2\sqrt{\log m}.
\]
Therefore, again by Chebyshev's inequality,
\[
\begin{aligned}
\Q_m(S_m>u'\log m+A\sqrt{\log m}) &\le \Q_m\left( \left|S_m-\E_{\Q_m}\left[S_m\right]\right| \ge u'\log m+A\sqrt{\log m}-\E_{\Q_m}\left[S_m\right] \right)\\
&\le \Q_m\left( \left|S_m-\E_{\Q_m}\left[S_m\right]\right| \ge \frac A2\sqrt{\log m} \right)\\
&\le \frac{\operatorname{Var}_{\Q_m}(S_m)} {\left(\frac A2\sqrt{\log m}\right)^2} = \frac{4\left(u'\log m+O_{u'}(1)\right)} {A^2\log m}\\
&= \frac{4u'}{A^2}+o(1).
\end{aligned}
\]
Taking \(A\) large enough, we get
\[
\liminf_{m\to\infty}\Q_m(A_m) \ge 1-\frac{4u'}{A^2}>0.
\]
Using the change of measure,
\[
\P(S_m\ge u\log m) \ge \P(A_m) = \E_{\Q_m}\left[ \frac{d\P}{d\Q_m}\mathbf 1_{A_m} \right] = \E\left[z'^{S_m}\right] \E_{\Q_m}\left[ z'^{-S_m}\mathbf 1_{A_m} \right].
\]
Since \(S_m\le u'\log m+A\sqrt{\log m}\) on \(A_m\) and \(z'>1\),
\[
\P(S_m\ge u\log m) \ge \E\left[z'^{S_m}\right] z'^{-(u'\log m+A\sqrt{\log m})} \Q_m(A_m).
\]
Taking logarithms and dividing by \(\log m\), we obtain
\begin{equation}
\label{right-E}
\liminf_{m\to\infty} \frac1{\log m} \log \P(S_m\ge u\log m) \ge 2(z'-1)-u'\log z'.
\end{equation}
Using \(z'=u'/2\), the right-hand side of~\eqref{right-E} becomes
\[
2\left(\frac{u'}2-1\right) - u'\log\frac{u'}2 = -\left(u'\log\frac{u'}2-u'+2\right) = -I(u').
\]
Thus, letting \(u'\downarrow u\), and using continuity of \(I\), we have
\[
\liminf_{m\to\infty} \frac1{\log m} \log \P(S_m\ge u\log m) \ge -I(u).
\]
Together with the upper bound, this proves
\[
\P(S_m\ge u\log m) = m^{-I(u)+o(1)}.
\]
\end{proof}

\subsection{Proof of Corollary~\ref{cor:directionless-local-finiteness}}
\label{app:proof-directionless-local-finiteness}

\begin{proof}
By Lemma~\ref{lem:directionless-finite-interval-incident-edges},
\[
N_{G_\infty^{(K)}}(v) \subseteq [v-Q_K(v),v+Q_K(v)].
\]
Since \(Q_K(v)<\infty\) almost surely, every fixed vertex has finite degree almost surely.
A countable intersection over \(v\in\Z\) completes the proof.
\end{proof}

\subsection{Proof of Lemma~\ref{lem:directionless-finite-interval-rooted-balls}}
\label{app:proof-directionless-finite-interval-rooted-balls}

\begin{proof}
For every \(x\in B_r(G_\infty^{(K)},v)\), we have
\[
[x-Q_K(x),x+Q_K(x)] \subseteq [v-\mathcal R_r(v),v+\mathcal R_r(v)] \subseteq[n].
\]
Lemma~\ref{lem:directionless-finite-interval-incident-edges} therefore gives
\[
N_{G_n^{(K)}}(x) = N_{G_\infty^{(K)}}(x), \qquad \forall x\in B_r(G_\infty^{(K)},v).
\]
Starting from the common root \(v\), induction on graph distance shows that the two graphs have the same vertices at distance at most \(r\) from \(v\).
Since the neighbor sets agree for every vertex in the radius-\(r\) ball, the induced edge sets agree as well.
Hence
\[
B_r(G_n^{(K)},v) = B_r(G_\infty^{(K)},v).
\]
\end{proof}

\section*{Acknowledgments}
This research was conducted at the UMN Duluth REU with support from Jane Street Capital, Axiom Math, University of Michigan, and a donation from Larry Penn. 
I am grateful to Joe Gallian and Colin Defant for running this program.
I would like to thank Mitchell Lee and Glenn Bruda for proofreading and providing insightful comments, and Michael Ren for his guidance during the early stages of this work. 
I would also like to thank my advisors, Eliot Hodges, Claire Kaneshiro, Noah Kravitz, and Carl Schildkraut, for helpful discussions and feedback.

\section*{AI Statement}
All proofs in this paper are the author's original work. Some ideas were inspired by interactions with LLMs, but the LLMs were unable to prove the main results.
The paper was completed in July 2026. The author used Claude Opus (versions 4.8 and 5) in preparing this manuscript, as follows, in chronological order:
(1) At the start of the project, the author used it to help read related papers.
(2) During the write-up stage, the author used it to write code for testing various candidate arguments and intuitions for Theorem~\ref{thm:max-degree-asymptotics} and Conjecture~\ref{conj:directionless-K-degree-tail}, and to search for references.
When developing proofs, the author also occasionally consulted Claude on the mathematical viability of candidate proof flows that the author thought might work, and adopted their suggestions selectively.
(3) After the draft was completed, non-generative AI tools (e.g. Overleaf AI assistant) were used to refine the writing.

\printbibliography
\end{document}